\documentclass[12pt,reqno,a4paper]{amsart}    

\usepackage[totalwidth=500pt, totalheight=710pt]{geometry}

\usepackage{amsmath, amssymb}
\usepackage{amsthm}

\usepackage{array}  
\usepackage{diagrams} 
\usepackage{graphicx}                      

\usepackage{bbm}         

\usepackage{mathrsfs}    
\renewcommand\mathcal\mathscr  

\usepackage{accents}     

\numberwithin{equation}{section}

\theoremstyle{plain}            
\newtheorem{theorem}{Theorem}[section]
\newtheorem*{maintheoremA*}{Main Theorem A} 
\newtheorem*{maintheoremB*}{Main Theorem B}
\newtheorem*{maintheoremC*}{Main Theorem C}
\newtheorem{proposition}[theorem]{Proposition}
\newtheorem{lemma}[theorem]{Lemma}

\theoremstyle{definition}       
\newtheorem{definition}[theorem]{Definition}

\newtheorem{example}[theorem]{Example}

\theoremstyle{remark}           
\newtheorem{remark}[theorem]{Remark}

\newcommand{\Sec}[1]{Section~\ref{sec:#1}}
\newcommand{\Secs}[2]{Sections~\ref{sec:#1} and~\ref{sec:#2}}

\newcommand{\Eq}[1]{Eq.~\eqref{eq:#1}}

\newcommand{\Fig}[1]{Figure~\ref{fig:#1}}

\newcommand{\Figs}[2]{Figures~\ref{fig:#1} and~\ref{fig:#2}}

\newcommand{\Thm}[1]{Theorem~\ref{thm:#1}}

\newcommand{\Ex}[1]{Example~\ref{ex:#1}}
\newcommand{\Exs}[2]{Examples~\ref{ex:#1} and~\ref{ex:#2}}
\newcommand{\ExS}[2]{Examples~\ref{ex:#1}--\ref{ex:#2}}

\newcommand{\Lem}[1]{Lemma~\ref{lem:#1}}

\newcommand{\Prp}[1]{Proposition~\ref{prp:#1}}
\newcommand{\Prps}[2]{Propositions~\ref{prp:#1} and~\ref{prp:#2}}
\newcommand{\Prpss}[3]{Propositions~\ref{prp:#1},~\ref{prp:#2}
  and~\ref{prp:#3}}

\newcommand{\Rem}[1]{Remark~\ref{rem:#1}}

\newcommand{\Def}[1]{Definition~\ref{def:#1}}

\newcommand{\dd}    {\, \mathrm d}    

\DeclareMathOperator{\dom}    {dom}

\DeclareMathOperator{\id}     {id}  

\DeclareMathOperator{\ran}    {ran}

\DeclareMathOperator{\supp}   {supp}

\newcommand{\specsymb} {\sigma} 

\newcommand{\spec}[2][{}]   {\specsymb_{\mathrm{#1}}(#2)}

\def\XXint#1#2#3{{\setbox0=\hbox{$#1{#2#3}{\int}$}
     \vcenter{\hbox{$#2#3$}}\kern-.5\wd0}}
\def\XXsum#1#2#3{{\setbox0=\hbox{$#1{#2#3}{\int}$}
     \vcenter{\hbox{$#2#3$}}\kern-.60\wd0}}

\newcommand{\R}{\mathbb{R}} 
\newcommand{\C}{\mathbb{C}} 
\newcommand{\N}{\mathbb{N}} 
\newcommand{\Z}{\mathbb{Z}} 

\renewcommand{\phi}{\varphi}   
\newcommand{\e}{\mathrm e}  
\newcommand{\im}{\mathrm i} 

\newcommand{\wt}{\widetilde}           
\newcommand {\qf}[1]{\mathfrak{#1}}    

\newcommand{\HS}{\mathcal H}           

\newcommand{\Sobsymb} {\mathsf H}      
\newcommand{\Contsymb} {\mathsf C}     
\newcommand{\Lsymb}    {\mathsf L}     
\newcommand{\lsymb}    {\ell}          

\newcommand{\Sobspace}[1]{\Sobsymb^{#1}}      
\newcommand{\Lsqrspace}    {\Lsymb_2}     
\newcommand{\lsqrspace}    {\lsymb_2}          

\newcommand{\Cont}[2][{}]{\Contsymb^{#1}({#2})}

\newcommand{\Lsqr}[2][{}]{\Lsymb_2^{#1}({#2})} 
 
\newcommand{\lsqr}[2][{}]{\lsymb_2^{#1}({#2})}   

\newcommand{\Sob}[2][1]{\Sobsymb^{#1}({#2})}         
\newcommand{\Sobn}[2][1]{\ring \Sobsymb^{#1}({#2})}  

\newcommand{\Sobx}[3][1]{\Sobsymb_{{#2}}^{#1}({#3})} 
  
\newcommand{\abs}[2][{}]{\lvert{#2}\rvert_{{#1}}}    
\newcommand{\abssqr}[2][{}]{\lvert{#2}\rvert^2_{#1}} 

\newcommand{\normsqr}[2][{}]{\|{#2}\|^2_{#1}} 

\newcommand{\iprod}[3][{}]{\langle{#2},{#3}\rangle_{#1}}  

\newcommand{\set}[2]{\{ \, #1 \, | \, #2 \, \} }      
\newcommand{\bigset}[2]{\bigl\{ \, #1 \, \bigl|\bigr. \, #2 \, \bigr\} }

\newcommand{\map}[3]{ #1 \colon #2 \longrightarrow #3}    

\newcommand{\bd}  {\partial}                
\newcommand{\clo}[1]{\overline{{#1}}} 

\newcommand{\dcup}{\mathbin{\mathaccent\cdot\cup}}

\DeclareMathOperator*{\bigdcup}{\mathaccent\cdot{\bigcup}}

\newcommand{\restr}[1]{{\restriction}_{#1}} 

\newcommand{\orth}{\bot}                    

\newcommand{\1}{\mathbbm 1}                    
\newcommand{\und}{\quad\text{and}\quad}
\newcommand{\Und}{\qquad\text{and}\qquad}

\newcommand{\Neu}{{\mathrm N}}              
\newcommand{\Dir}{{\mathrm D}}              
\newcommand{\laplacian}[1][{}]{\Delta^{#1}} 
\newcommand{\dlaplacian}[1][{}]{\check\Delta^{{#1}}}
\newcommand{\discr}[1]{\check {#1}}         

\newcommand{\qlaplacian}[1][{}]{\Delta^{#1}} 

\newcommand{\orient}[1]{\accentset{\curvearrowright}{#1}} 

\newcommand{\mc}{\mathcal}

\newcommand{\wh}{\widehat}

\newcommand{\color}[2][{}]{}        

\newcommand{\rel}[1]{#1_{\mathrm r}}   
\newcommand{\unor}[1]{\bar #1}         

\begin{document}

\title[Eigenvalue bracketing for graphs]{Eigenvalue bracketing for
  discrete and metric graphs}

\author{Fernando Lled\'o}
 \address{Department of Mathematics, University Carlos~III Madrid,
     Avda.~de la Universidad~30, E-28911 Legan\'es (Madrid),
     Spain and Institut f\"ur Reine und Angewandte Mathematik,
        RWTH-Aachen University,
        Templergraben 55,
        D-52062 Aachen,
        Germany (on leave)}
   \email{flledo@math.uc3m.es \emph{and} lledo@iram.rwth-aachen.de}

\author{Olaf Post} 
\address{Institut f\"ur Mathematik,
         Humboldt-Universit\"at, Rudower Chaussee~25,
         12489 Berlin, Germany}
         \email{post@math.hu-berlin.de}

\date{Compiled on \today}


\begin{abstract}
  We develop eigenvalue estimates for the Laplacians on discrete and
  metric graphs using different types of boundary conditions at the
  vertices of the metric graph. Via an explicit correspondence of the
  equilateral metric and discrete graph spectrum (also in the
  ``exceptional'' values of the metric graph corresponding to the
  Dirichlet spectrum) we carry over these estimates from the metric
  graph Laplacian to the discrete case.  We apply the results to
  covering graphs and present examples where the covering graph
  Laplacians have spectral gaps.
 \end{abstract}

\maketitle


%

\section{Introduction}
\label{sec:intro}
%

Analysis on graphs is an area of current research in mathematics with many
applications e.g.\ in network theory, nano-technology, optics, chemistry and
medicine. In this context one studies different kinds of linear operators,
typically Laplacians, on a graphs. From the spectral properties of these
operators one may infer relevant information of the corresponding model. For
example, the tight binding model in physics describes atoms and molecules by a
nearest neighbour model closely related to the discrete graph Laplacian.
Moreover, network properties like connectivity can be described with spectral
graph theory.  In applications, the spectrum may encode transport properties
of the medium.  We will call an interval \emph{disjoint} from the spectrum a
\emph{spectral gap}. In applications, a spectral gap may describe a set of
wave-lengths for which no transport is permitted through the media.

There are basically two ways to give a ``natural'' definition of the
Laplace operator on graphs: first, on discrete graphs, the operator
acts on functions on the \emph{vertices} as \emph{difference}
operator. Here, edges play a secondary role as labels that connect the
vertices.  Second, one can consider the graph as a (non-discrete)
metric space consisting of vertices and edges as
\emph{one-dimensional} spaces.  In this context one defines
\emph{differential} operators acting on functions on the edges.
Laplacians are second order operators with suitable boundary
conditions on the vertices chosen in such a way that the operator is
self-adjoint in the corresponding $\Lsqrspace$-space.  One usually
refers to a metric graph together with a self-adjoint differential
operator as a \emph{quantum graph}.  Recent interesting reviews on
discrete geometric analysis and quantum graphs can be found
in~\cite{sunada:pre08} resp.~\cite{kuchment:pre08} (see also
references therein).

The aim of the present paper is to use spectral results for the metric
graph to obtain spectral information of the discrete Laplacian. In
particular, we will obtain results on the spectrum of infinite
discrete covering graphs.  This partially answers a question of Sunada
concerning the spectrum of infinite discrete graphs
\cite[p.~64]{sunada:07}.  In particular, we generalise the so-called
\emph{Neumann-Dirichlet bracketing} (see below) to the Laplacian
acting on a metric graph, where the lower bound estimate of the
Neumann eigenvalue is replaced by the Kirchhoff condition. Due to an
explicit relation between the spectrum of the discrete and
(equilateral) metric Laplacian, we can carry over the eigenvalue
estimates to the discrete case. We also treat the exceptional
eigenvalues in this relation (usually due to the Dirichlet spectrum of
a single edge), and relate them with relative homology of the graph
and its boundary.  This gives a complete relation between the discrete
and metric spectra (see Theorem~A below).

\subsection*{The basic idea of the eigenvalue bracketing}
Our basic technique is to localise the eigenvalues within suitable
closed intervals which we can control. We call this process
\emph{bracketing}. Since this technique is crucial for our analysis,
we will briefly recall the main idea here.  

Dirichlet-Neumann bracketing is a tool usually available for
differential operators like Schr\"odinger operators or Laplacians on
manifolds. The simplest example is provided by the operator
$\laplacian f=-f''$ on the interval $[0,1]$. In order to obtain a
self-adjoint operator in $\Lsqr{0,1}$ one has to fix boundary
conditions at $0$ and $1$.  
A very elegant way to provide such
conditions is to define the Laplacian via an associated quadratic
form
\begin{equation*}
  \qf h(f) := \int_0^1 \abssqr {f'(x)} \dd x, \quad f \in \dom \qf h
  \qquad \text{related by} \quad \iprod f{\laplacian f} = \qf h(f).
\end{equation*}
The quadratic form domain is a closed subspace of the Sobolev space
$\Sob{0,1}$. The two extremal cases are
\begin{enumerate}
\item the \emph{Dirichlet} boundary condition, $\dom \qf h^\Dir := \set
  {f \in \Sob{0,1}} {f(0)=0, \; f(1)=0}$.
\item the \emph{Neumann} boundary condition, $\dom \qf h^\Neu :=
  \Sob{0,1}$.
\end{enumerate}
Note that the usual Neumann conditions $f'(0)=0$ and $f'(1)=0$ only
enter in the \emph{operator} domain by requiring the boundary terms to
vanish which appear after partial integration. For details, we refer
to~\cite[Sec.~VIII.6]{reed-simon-1}
and~\cite[Sec.~XIII.15]{reed-simon-4} or~\cite{davies:95}.  Any other
(linear) boundary condition, like e.g.\ the $\vartheta$-equivariant
condition $f(1)=\e^{\im \vartheta} f(0)$ leads to a space $\dom \qf
h^\vartheta$ between $\dom \qf h^\Dir$ and $\dom \qf h^\Neu$ (the
action of $\qf h^\vartheta$ being the same, namely $\qf
h^\vartheta(f)= \normsqr{f'}$). Floquet theory implies that the
spectrum of the corresponding ($\Z$-periodic) Laplacian
$\laplacian_\R$ on $\R$ is given by $\set{\lambda_k^\vartheta}{k \in
  \N, \vartheta \in [0,2\pi]}=[0,\infty)$.  The variational
characterisation of the associated eigenvalues is given by
\begin{equation}
  \label{eq:min.max}
  \lambda_k^\bullet = 
     \inf_{D \subset \dom \qf h^\bullet}
     \sup_{f \in D} \frac{\qf h^\bullet(f)} {\normsqr f}
\end{equation}
where $D$ runs through all $k$-dimensional subspaces and the dot
$^\bullet$ is a placeholder for the labels $\mathrm N$, $\mathrm D$,
$\vartheta$.  Extending the non-negative forms $\qf h^\bullet$
naturally to the whole Hilbert space by $\qf h^\bullet(f):= \infty$ if
$f \notin \dom \qf h^\bullet$, the extended forms become monotone in
the obvious sense, i.e.  $\qf h^\Neu(f) \le \qf h^\vartheta(f) \le \qf
h^\Dir(f)$ for all $f \in \Lsqr {0,1}$ (opposite to the inclusion of
the domains). It follows now from \Eq{min.max} that
\begin{equation*}
  \lambda_k^\Neu \le \lambda_k^\vartheta \le \lambda^\Dir,
\end{equation*}
to what we will refer to as \emph{bracketing}.  In this simple example
the bracketing does not imply the existence of spectral gaps of
$\laplacian_\R$ inside $[0,\infty)$, since
$\lambda_k^\Neu=\pi^2(k-1)^2$ and $\lambda_k^\Dir=\pi^2k^2$, $k=1,2,
\dots$, and therefore the intervals $I_k:= [\lambda_k^\Neu,
\lambda_k^\Dir]$ cover already $[0,\infty)$. Of course, we do not
expect gaps here since $\spec {\laplacian_\R}=[0,\infty)$.

The strength of this bracketing method can be seen in \Prp{sp.incl}
where we use the same idea for eigenvalues of equilateral metric graph
Laplacians and \emph{arbitrary} finite-dimensional unitary
representations $\rho$.  \Prp{sp.incl} may be seen as the core of our
analysis. Its proof is amazingly simple, namely, it is a vector-valued
generalisation of the above argument.

\subsection*{Main results}

Let us briefly describe our main results: Denote by $\mc N(\lambda)$
the eigenspace of the standard (Kirchhoff) metric Laplacian, and by
$\discr{\mc N}(\mu)$ the eigenspace of the standard discrete Laplacian
(for precise definitions, see \Secs{dg}{mg}).

The following theorem resumes \Prpss{sp.rel}{dec.spaces}{hom.ef},
where the precise statements can be found. The first statement for
equilateral graphs (i.e., metric graphs with constant length function,
say, $\ell_e=1$) is standard (see
e.g.~\cite{von-below:85,exner:97b,cattaneo:97,pankrashkin:06a,post:pre07c,
  bgp:08}) and only mentioned for completeness:
\begin{maintheoremA*}[\Prpss{sp.rel}{dec.spaces}{hom.ef}]
  Let $X$ be a compact, connected and equilateral metric graph and set
  $\mu(\lambda):=1-\cos(\sqrt \lambda)$.
  \begin{enumerate}
  \item If $\lambda \notin \set{n^2\pi^2}{n=1,2,\dots}$, then there is
    an isomorphism $\map {\Phi_\lambda} {\discr{\mc N}(\mu(\lambda))}
    {\mc N (\lambda)}$. The corresponding metric eigenfunctions are
    called \emph{vertex based}.
  \item If $\lambda_n = n^2 \pi^2$ and $n$ even, then there is an
    injective homomorphism $\map {\Psi_n} {H_1(X)} {{\mc
        N}(\lambda_n)}$, where $H_1(X)$ is the first homology group.
    The range of $\Psi_n$ consists of functions vanishing on all
    vertices (``Dirichlet eigenfunctions''), called \emph{edge-based}
    or \emph{topological} eigenfunctions of the metric graph.  The
    orthogonal complement of the range of $\Psi_n$ contains an
    additional eigenfunction $\phi_n$ which is constant as function
    restricted to the set of vertices, called \emph{trivial vertex
      based}.
  \end{enumerate}
\end{maintheoremA*}
For shortness, we omit the case $n$ odd, in which a similar statement
with $H_1(X)$ replaced by the ``unoriented'' homology group $\unor
H_1(X)$ holds. In this case, one has to distinguish whether $G$ is
bipartite or not.  In the former case, the orthogonal complement of
the range of $\Psi_n$ contains the additional eigenfunction $\phi_n$
related to the discrete bipartite eigenfunction. In the latter case,
$\Psi_n$ is already an isomorphism.  Moreover, a similar result holds
when we consider Laplacians with Dirichlet conditions on a subset $\bd
V$ of the vertices. In this case, the \emph{relative} homology group
$H^1(X,\bd V)$ enters.  The multiplicities of the eigenvalues were
already calculated in~\cite{von-below:85} by a direct proof without
using the homology groups. The advantage of using homology groups is
that is can be generalised to other types of vertex boundary
conditions (like Dirichlet and equivariant) in a natural way, see
e.g.\ \Rem{twisted.hom}).

Let now $X \to X_0$ be a covering of metric graphs (i.e., a covering
respecting the combinatorial graph structure and the length function).
For the next statement, the metric graph need not
to be equilateral.
\begin{maintheoremB*}[\Thm{brack.mg}]
  Let $X \to X_0$ be a covering of metric graphs with compact quotient and
  residually finite covering group $\Gamma$ and denote by $\qlaplacian_X$ the
  Kirchhoff Laplacian.  Then
  \begin{equation*}
    \spec{\qlaplacian_X} \subset \bigcup_{k \in \N} I_k, \qquad
      I_k=\bigl[\lambda_k,\lambda_k^{\bd V}\bigr],
  \end{equation*}
  where $\lambda_k$ and $\lambda_k^{\bd V}$ are the eigenvalues of the
  Kirchhoff and Kirchhoff-Dirichlet Laplacian on a fundamental domain $Y
  \subset X$.  In particular, for any subset $M \subset [0,\infty)$ such that
  $M\cap \bigcup_k I_k = \emptyset$, then
  $M\cap\spec{\qlaplacian_X}=\emptyset$.
\end{maintheoremB*}
Abelian groups, finite extensions of Abelian groups (so-called
\emph{type-I-groups}) and free groups are examples of the large class of
residually finite groups (see~\cite{lledo-post:08} for more details). For
Abelian groups, the Floquet-Bloch decomposition can be used in order to
calculate the spectrum of the operator on the covering, leading to a detailed
analysis in certain models, see e.g.~\cite{kuchment-post:07} for hexagonal
lattices (modeling carbon nano-structures).

We refer to the intervals $I_k=I_k(Y,\bd V)$ as
\emph{Kirchhoff-Dirichlet (KD) intervals}. Note that they depend
usually on the fundamental domain. The Kirchhoff condition plays the
role of the Neumann condition in the usual \emph{Dirichlet-Neumann}
bracketing. Note that the Kirchhoff condition is optimal in a sense
made precise in \Rem{kir.opt}, namely that a symmetrised version of
the KD intervals (explained below) give the exact spectrum of the
corresponding (Abelian) covering Laplacian.

We call the set $M$ also a \emph{spectral gap}.  Note that we do not
assume that the spectral gap is \emph{maximal}, i.e., if we state that
the spectrum has two disjoint gaps $(a_1,b_1)$ and $(a_2,b_2)$ with
$b_1 \le a_2$ we do not make a statement about the existence of
spectrum inside $[b_1,a_2]$. In certain situations (e.g.\ if $\Gamma$
is amenable), we can assure the existence of spectrum between the
gaps, and therefore have a lower bound on the number of components of
$\spec {\qlaplacian_X}$ in terms of the components of $\bigcup_k I_k$
(see \Thm{amenable}).

For an equilateral metric graph, we can combine the last two theorems
and obtain the following \emph{discrete Kirchhoff-Dirichlet
  bracketing}. Let $G \to G_0$ be a covering of discrete graphs with
fundamental domain $H$ (being a subgraph of $G$ with vertex set $V(H)$
and boundary $\bd V$):
\begin{maintheoremC*}[\Thm{brack.dg}]
  Assume that the covering group is residually finite, then
  \begin{equation*}
    \spec{\dlaplacian_G} \subset \bigcup_{k=1}^{\abs {V(H)}} J_k, \qquad
      J_k=[\mu_k,\mu_k^{\bd V}],
  \end{equation*}
  where $\mu_k$ and $\mu_k^{\bd V}$ are the eigenvalues of the
  discrete Laplacians on the fundamental domain $H$ with Dirichlet
  condition on $\bd V$ in the latter case.  In particular, for any
  subset $M \subset [0,\infty)$ such that $M \cap \bigcup_k J_k =
  \emptyset$, then $M\cap\spec{\dlaplacian_X}=\emptyset$.
\end{maintheoremC*}
We refer to the intervals $J_k=J_k(H,\bd V)$ as the \emph{discrete
  Kirchhoff-Dirichlet intervals}. Note that this method allows to
determine in a very easy way whether a set $M$ is not contained in the
spectrum of the covering Laplacian. The only step to be done is to
calculate the eigenvalues $\mu_k$ and $\mu_k^{\bd V}$ (which give
immediately the corresponding metric eigenvalues for equilateral
graphs) and check whether neighbouring KD intervals $J_k$ have empty
intersection. We will see in \Sec{ex}, that in simple examples, only
the first KD intervals do not overlap.  As in the case of manifolds
and Schr\"odinger operators (see
e.g.~\cite{hempel-post:03,lledo-post:07,lledo-post:08}) we expect that
the number of gaps should be large if the fundamental domain has
``small'' boundary $\bd V$ compared to the number of vertices $V(H)$
and edges $E(H)$ inside. In other words, a ``high contrast'' between
the different copies of a suitable fundamental domain is necessary in
order that our method works.

It is a priori not clear how the eigenvalue bracketing can be seen
directly for discrete Laplacians, so our analysis may serve as an
example of how to use metric graphs to obtain results for discrete
graphs.

\subsection*{Structure of the article}
The structure of the paper is as follows: in the following two
sections we present the basic definitions and results for various
Laplacians on discrete and metric graphs.  \Secs{sp.rel}{hom} contains
the complete relation of discrete and equilateral metric graphs and in
particular Main Theorem~A. Details on the different homologies needed
can be found in \Sec{hom}.  \Sec{ev.brack} is devoted to the
definition of the metric and discrete Kirchhoff-Dirichlet intervals
and contains a careful analysis of the metric eigenvalues including
multiplicities.  \Sec{eq.lapl} contains relevant information on
equivariant Laplacians and the basic idea of decoupling an equivariant
Laplacian via Dirichlet and Kirchhoff Laplacians
(see~\Prps{sp.incl}{sp.incl.dg}).  In \Sec{cov.gr} we combine the
results on equilateral Laplacians and KD intervals in order to prove
our Main Theorems~B and~C. The last section provides several examples
of graphs with spectral gaps.

%
\section{Discrete graphs}
\label{sec:dg}
%
Let $G=(V,E,\bd)$ be a (connected) discrete graph, i.e., $V=V(G)$ is
the set of vertices, $E=E(G)$ the set of edges and $\map {\bd=\bd_G} E
{V\times V}$ the connection map, $\bd e = (\bd_-e,\bd_+e)$ is the pair
of the initial and terminal vertex, respectively. Clearly, $\bd_\pm e$
fixes an orientation of the edge $e$.  We prefer to consider $E$ and
$V$ as independent sets (and not the edge sets as pairs of vertices),
in order to treat easily multiple edges (i.e., edges $e_1$, $e_2$ with
$\{\bd_-e_1,\bd_+e_1\}=\{\bd_-e_2,\bd_+e_2\}$) and self-loops (i.e.,
edges with $\bd_- e = \bd_+e$).  For two subsets $A, B \subset V$ we
denote by
\begin{equation*}
  E^+(A,B) := \set {e  \in E}{\bd_-e \in A,\, \bd_+e \in B }
\end{equation*}
the set of edges with terminal vertex in $A$ and initial vertex in
$B$, and similarly, $E^-(A,B):= E^+(B,A)$. Moreover we let $E(A,B):=
E^+(A,B) \dcup E^-(A,B)$ be the \emph{disjoint} union of all edges
between $A$ and $B$. Due to the disjoint union, a self-loop at a
vertex $v \in A \cap B$ is counted twice in $E(A,B)$.  In particular,
$E(v,w)$ is the set of all edges between the vertices $v$ and $w$; and
\begin{equation*}
  E_v^\pm := E^\pm(V,v) =\set{e \in  E}{\bd_\pm e=v}
\end{equation*}
is the set of edges terminating ($+$) and starting ($-$) at $v$.
Similarly, $E_v = E_v^+ \dcup E_v^-$ is the set of all edges at $v$.
We call
\begin{equation*}
  \deg v := |E_v|
\end{equation*}
the \emph{degree} of the vertex $v$ in the graph $G$. Note that a
self-loop at the vertex $v$ increases the degree by $2$.

A graph is called \emph{bipartite} if there is a disjoint
decomposition $V=A \dcup B$ such that $E=E(A,B)$, i.e., if each edge
has exactly one end-point in $A$ and the other in $B$.

We will use frequently the following elementary fact about reordering
a sum over edges and vertices, namely
\begin{equation}
  \label{eq:reorder}
  \sum_{e \in E} F(\bd_\pm e,e) 
   = \sum_{v \in V} \sum_{e \in E_v^\pm} F(v,e)
\end{equation}
for a function $(v,e) \mapsto F(v,e)$ depending on $v$ and $e \in E_v$
with the convention that a sum over the empty set is $0$. Note that
this equation is also valid for self-loops and multiple edges. The
reordering is a bijection since the union $E = \bigdcup_{v \in V}
E_v^\pm$ is \emph{disjoint}.

We start with a more general setting, namely with \emph{weighted}
graphs, i.e., we assume that there are two functions $\map{m=m_V} V
{(0,\infty)}$ and $\map{m=m_E}E{(0,\infty)}$ (mostly denoted by the
same symbol $m$) associating to a vertex $v$ its weight $m(v)$ and to
an edge $e$ its weight $m_e$. We will call $(G,m)$ a \emph{weighted}
discrete graph. The basic Hilbert spaces associated with $(G,m)$ are
\begin{align*}
  \lsqr {V,m} &:= \bigset{ \map F V \C}
  { \normsqr[V,m] F := \sum_{v \in V} |F(v)|^2 m(v) < \infty},\\
  \lsqr {E,m} &:= \bigset{ \map \eta V \C}
  { \normsqr[V,m] f := \sum_{e \in E} |\eta_e|^2 m_e < \infty}.
\end{align*}
We define the \emph{discrete exterior derivative} $d$ as
\begin{equation*}
  \map d {\lsqr{V,m}}{\lsqr{E,m}}, \qquad
  (dF)_e = F(\bd_+ e) - F(\bd_-e).
\end{equation*}
We define the \emph{relative} weight $\map{\rho} V {(0,\infty)}$ as
\begin{equation}
  \label{eq:rel.weight}
  \rho(v) := \frac 1 {m(v)} \sum_{e \in E_v} m_e
\end{equation}
and we will assume throughout this article that
\begin{equation}
  \label{eq:rel.bdd}
  \rho_\infty := \sup_{v \in V} \rho(v) < \infty,
\end{equation}
i.e., that the relative weight is uniformly bounded. We will call the
weights \emph{normalised} if $\rho(v)=1$ for all vertices.  A
straightforward calculation using~\eqref{eq:reorder} shows that $d$ is
an operator with norm bounded by $(2 \rho_\infty)^{1/2}$. Similarly,
one can calculate the adjoint $\map{d^*}{\lsqr{E,m}}{\lsqr{V,m}}$ and
one gets
\begin{equation*}
  (d^*F)(v) = \frac 1 {m(v)} \sum_{e \in E_v} m_e \orient \eta_e(v),
\end{equation*}
where
\begin{align}
  \label{eq:def.orient}
  \orient \eta_e(v) &= \eta_e \quad \text{if $v=\bd_+ e$} &\text{and}&&
  \orient \eta_e(v) &=  -\eta_e \quad \text{if $v=\bd_- e$.}
\end{align}

The \emph{discrete Laplacian} is now defined as
\begin{equation}
  \label{eq:def.disc.lapl}
  \map{\dlaplacian = \dlaplacian_{(G,m)} := d^* d} {\lsqr {V,m}} {\lsqr {V,m}}
\end{equation}
and acts as
\begin{equation}
  \label{eq:disc.lapl}
  (\dlaplacian_{(G,m)} F)(v) 
  = \rho(v) F(v) - \frac 1 {m(v)} \sum_{e \in E_v} m_e F(v_e),
\end{equation}
where $v_e$ denotes the vertex on the edge $e \in E_v$ opposite to
$v$. If no confusion arises we also denote the Laplacian simply by
$\dlaplacian$. The \emph{standard discrete Laplacian} is the Laplacian
associated with the weights $m(v)=\deg v$ and $m_e=1$.  We will often
refer to the standard weighted graph as $(G,\deg)$ or simply as $G$.
Note that these weights are normalised, i.e., that $\rho(v)=1$.

\begin{remark}
  \label{rem:orient}
  Note that as second order difference operator, the Laplacian does
  not see the orientation of the graph, whereas the discrete exterior
  derivative as first order operator depend on the orientation. We
  will define below an \emph{unoriented} version of the exterior
  derivative $\unor d$ that does not see the orientation. The
  corresponding (co-)homologies for $d$ and $\unor d$ will be useful
  in order to analyse exceptional metric graph eigenfunctions composed
  of antisymmetric and symmetric Dirichlet eigenfunctions on a single
  edge (see \Sec{hom}).
\end{remark}

A graph without multiple edges (i.e., $\abs{E(v,w)}\le 1$ for all $v,w
\in E$ is called \emph{simple}. In particular, $\bd$ is injective and
we can consider $E$ as a subset of $V \times V$. In this case, we also
write $v \sim w$ if $v$ and $w$ are connected by an edge.

One reason for considering graphs with arbitrary weights is the fact
that one can express the standard Laplacian on a graph with multiple
edges and self-loops equivalently by a Laplacian on a simple graph by
changing the weights. We will use multiple edges and self-loops in
\ExS{mult-edge}{self-loops} in order to generate gaps. Note that the
corresponding discrete exterior derivatives will of course differ, as
well as the topology of the graph.  Nevertheless, the reduction to
simple graphs is more convenient when calculating the spectrum of the
Laplacian.

\subsection{Multiple edges}

Assume that $G$ is a graph with the standard weights $m(v)=\deg v$,
$m_e=1$ and that $G$ has multiple edges. We can pass to a graph $\wt
G$ having the same set of vertices as $G$ but only simple edges.  The
multiple edges $e \in E(v,w)$ in $G$ are replaced by a single edge
$(v,w)$ (not taking care about the original orientation) in $\wt G$.
Note that for the degree $\deg_{\wt G} v \le \deg_G v$ where
$\deg_{\wt G} v$ denotes the degree of $v$ in the simple graph $\wt
G$. We define
\begin{equation*}
  \wt m(v) := \deg_G v   \Und   \wt m_{(v,w)} := \abs{E(v,w)},
\end{equation*}
where $\deg_G v$ is the degree in the original graph. Now, the relative
weight $\wt \rho$ is still normalised, since
\begin{equation*}
  \wt \rho(v) 
  = \frac 1 {\wt m(v)} \sum_{w \sim v} \wt m_{(v,w)}
  = \frac 1 {\deg_G v} \sum_{w \sim v} \abs{E(v,w)} = 1.
\end{equation*}
Note that the Laplacians on $(\wt G, \wt m)$ and $(G,\deg)$ agree.

\subsection{Self-loops}
Assume that $G$ is a graph with a self-loop $e$, i.e.,
$\bd_+e=\bd_-e=v$. Obviously, for such an edge, we have $(dF)_e=0$,
i.e., we can eliminate this edge from $E$. We define a new graph $\wt
G$ having again the same vertex set as $G$ and where the edge set $\wt
E$ is the original edge set without self-loops. The degree in the new
graph is given by $\deg_{\wt G} v = \deg_G v - \abs{E(v,v)}$, i.e.,
the original degree minus \emph{twice} the number of self-loops
removed (remember that $E(v,v)$ was defined as the formal disjoint
union of $E^+(v,v)$ and $E^-(v,v)$). We set
\begin{equation*}
  \wt m(v):= \deg_G v \Und   \wt m_e=1,
\end{equation*}
so that the relative weight $\wt \rho$ satisfies
\begin{equation*}
  \wt \rho(v) 
  = \frac 1 {\wt m(v)} \sum_{e \in \wt E_v} 1
  = \frac {\deg_G v - \abs{E(v,v)}} {\deg_G v} < 1
\end{equation*}
provided there was a self-loop at $v$. Again, the corresponding
Laplacians on $(G,\deg)$ and $(\wt G,\wt m)$ agree.

\subsection{Matrix representation of the Laplacian}
For concrete computations of the eigenvalues of the weighted
Laplacian, it is convenient to have the associated matrix at hand. Let
$\{\phi_v\}_v$ be the standard orthonormal basis of $\lsqr{V,m}$,
where $\phi_v(w):=m(v)^{-1/2}$ if $v=w$ and $\phi_v(w)=0$ otherwise.
Then the matrix $L$ associated to the Laplacian
$\dlaplacian=\dlaplacian_{(G,m)}$ is given as
\begin{equation}
  \label{eq:matrix}
  \dlaplacian_{v,w} := \iprod{\phi_v}{\dlaplacian \phi_w}
  =
  \begin{cases}
    \rho(v) - \dfrac 1 {m(v)}{\sum_{e \in E(v,v)}m_e} &
           \text{if $v=w$}\\
    - \dfrac 1 {(m(v)m(w))^{1/2}} {\sum_{e \in E(v,w)}m_e} &
           \text{if $v \sim w$, $v \ne w$,}\\
    0   & \text{otherwise.}
  \end{cases}
\end{equation}
If the graph has the standard weights, then we obtain
\begin{equation}
  \label{eq:mat.std}
  \dlaplacian_{v,w} := \iprod{\phi_v}{\dlaplacian \phi_w}
  =
  \begin{cases}
    \dfrac {\deg (v) - \abs{E(v,v)}} {\deg v} &
           \text{if $v=w$}\\
    -\dfrac {\abs{E(v,w)}} {(\deg v \deg w)^{1/2}} &
           \text{if $v \sim w$, $v \ne w$,}\\
    0   & \text{otherwise.}
  \end{cases}
\end{equation}
Note that the latter expression also applies for graphs with multiple
edges and loops, inserting as degree function the degree of the
original (non-simple) graph.

\subsection{Discrete Dirichlet Laplacians}

A \emph{boundary} of $G$ is a subset $\bd V$ of $V$. We denote by
$\ring V:= V \setminus \bd V$ its complement, the \emph{inner}
vertices.  We set
\begin{equation*}
    \lsqr[\bd V] {V, m} := \bigset{F \in \lsqr {V,m}}
  { F \restr {\bd V} = 0}
\end{equation*}
and define the \emph{Dirichlet} discrete exterior derivative $d_0$ as
the restriction of $d$ to $\lsqr[\bd V]{V,m}$. Formally, we can write
$d_0 := d \circ \iota$, where $\iota$ is the canonical embedding of
$\lsqr[\bd V]{V,m}$ into $\lsqr{V,m}$. The adjoint of $d_0$ is
$d_0^*=\iota^* \circ d^*$, i.e.,
\begin{equation*}
   d_0^* \eta = (d^*\eta) \restr{\ring V},
\end{equation*}
since $\iota^* F$ is the restriction of $F$ onto the inner vertices
$\ring V$.

The \emph{discrete Dirichlet Laplacian} is defined as
\begin{equation*}
  \dlaplacian[\bd V] = \dlaplacian[\bd V]_{(G,m)} := d_0^* d_0
\end{equation*}
and acts as in~\eqref{eq:disc.lapl}, but only for $v \in \ring V$.

\begin{remark}
  One can give an equivalent definition of the Dirichlet Laplacian as
  a discrete Laplacian on the graph $\ring G$ with vertex set $\ring
  V$ and edge set $\ring E := E \setminus E(V,\bd V)$ (removing the
  edges to the boundary or inside the boundary). Again, this leads to
  a weighted Laplacian: If for instance, $\dlaplacian[\bd V]_G$ is the
  Dirichlet Laplacian with standard weights, we define
  \begin{equation*}
    \ring m(v) := \deg_G v    \Und   \ring m_e := 1
  \end{equation*}
  having again a non-normalised relative weight $\rho(v) < 1$ provided
  $v$ is joined with $\bd V$ by an edge in the original graph $G$.
\end{remark}

\subsection{Bipartiteness and the spectrum}
Let us recall the following spectral characterisation of bipartiteness
of a graph:
\begin{proposition}
  \label{prp:bipartite}
  Let $(G,m)$ be a weighted, connected graph with normalised weights
  (i.e., $\rho=1$). Assume in addition that $G$ has finite mass
  $m(V)=\sum_{v \in V} m(v) < \infty$ (e.g.\ that $G$ is finite). Then
  the following assertions are equivalent:
  \begin{enumerate}
  \item
    \label{bip1}
    The graph $G$ is bipartite
  \item
    \label{bip2}
    If $\mu \in \spec{\dlaplacian_{(G,m)}}$ then $2-\mu \in
    \spec{\dlaplacian_{(G,m)}}$. For short, we write
    $\spec{\dlaplacian_{(G,m)}}=2-\spec{\dlaplacian_{(G,m)}}$. The
    multiplicity is preserved.
  \item $2$ is an eigenvalue of $\dlaplacian_{(G,m)}$.
  \end{enumerate}
  Moreover, the implication~$\eqref{bip1} \Rightarrow \eqref{bip2}$
  holds for the discrete \emph{Dirichlet} Laplacian $\dlaplacian[\bd
  V]_{(G,m)}$.  Similarly, if $m(V)$ is infinite, then the
  implication~$\eqref{bip1} \Rightarrow \eqref{bip2}$ is still valid.

\end{proposition}
\begin{proof}
  The proof of the equivalence for $\bd V=\emptyset$ and finite graphs
  can be found e.g.\ in~\cite{chung:97}; the case $\bd V \ne
  \emptyset$ follows similarly. If $G$ has finite mass, then the
  constant function $\1_V$ is in $\lsqr {V,m}$, and the argument for
  finite graphs carries over.

  If $m(V)$ is infinite, then the spectral symmetry follows from the
  fact that
  \begin{equation}
    \label{eq:bip}
    \dlaplacian T = T (2 - \dlaplacian), \qquad T := \1_A - \1_B
  \end{equation}
  where $V=A \dcup B$ is the bipartite partition. Here, $T$ is a
  unitary involution (i.e., $T=T^*=T^{-1}=T^2$) on $\lsqr {V,m}$.
\end{proof}
Note that in the finite mass case, $T$ interchanges the constant
eigenfunction and the eigenfunction $\1_A - \1_B$ associated to the
eigenvalue $2$, also called the \emph{bipartite eigenfunction}.
Moreover, the condition~\eqref{eq:bip} is equivalent to the fact that
$T$ \emph{anticommutes} with the \emph{principal part} of the
Laplacian $L:=\id - \dlaplacian$, i.e., that $\{L,T\}=LT+TL=0$.

\subsection{Unoriented exterior derivatives}
\label{sec:unor.ext.der}
We briefly describe another sort of discrete exterior derivative, this
time an operator which does not see the orientation of the graph. More
precisely, we define the \emph{unoriented discrete exterior
  derivative} as
\begin{equation*}
  \map {\unor d} {\lsqr{V,m}}{\lsqr{E,m}}, \qquad
  (dF)_e = F(\bd_+ e) + F(\bd_-e),
\end{equation*}
i.e., compared with the (oriented) version $d$, we only change the
sign of the value of $F$ at the initial vertex. As a consequence, the
corresponding adjoint is given by
\begin{equation*}
  (\unor d^*F)(v) = \frac 1 {m(v)} \sum_{e \in E_v} m_e \eta_e.
\end{equation*}
One can also define a Laplacian associated via $\unor \dlaplacian :=
\unor d^* \unor d$, and the relation with the Laplacian $\dlaplacian =
d^* d$ is given by
\begin{equation}
  \label{eq:unor.lapl}
  \unor \dlaplacian = 2 \rho - \dlaplacian,
\end{equation}
where $\rho$ denotes the multiplication operator with the relative
weight. We will need the operators $\unor d$ and $\unor d^*$ in
\Sec{hom}. For more details and a general concept, in which the
oriented and unoriented version of an exterior derivative embed
naturally, we refer to~\cite{post:pre07a} (see
also~\cite{post:pre07c,post:pre07d}).

As for the oriented exterior derivative, we can also define a
Dirichlet version of $\unor d$, namely,
\begin{equation*}
  \map{\unor d_0} {\lsqr[\bd V]{V,m}} {\lsqr {E,m}}, \qquad
  \unor d_0 := \unor d \circ \iota.
\end{equation*}
As before, we have $\unor d_0^* \eta = (\unor d^*\eta) \restr{\ring V}$
for the adjoint.

%
\section{Metric graphs}
\label{sec:mg}
%
Let $G=(V,E,\bd)$ be a discrete graph. A \emph{topological graph}
associated to $G$ is a CW complex $X$ containing only $0$-cells and
$1$-cells. The $0$-cells are the vertices $V$ and the
$1$-cells are labelled by the edge set $E$.

A \emph{length function} $\map \ell E {(0,\infty)}$ of a graph $G$ is
the inverse of an edge weight function $m$, i.e., $\ell_e=1/m_e$.  We
will assume that the edge weight is bounded, i.e., that there exists
$\ell_0 > 0$ such that
\begin{equation}
  \label{eq:len.bd}
  \ell_e \ge \ell_0, \qquad \forall \; e \in E.
\end{equation}

The \emph{metric graph} $X$ associated to a weighted discrete graph
$(G,m)$ is a topological graph associated to $(V,E,\bd)$ such that for
every edge $e \in E$ there is a continuous map $\map{\Phi_e}{\clo
  I_e}X$, $I_e:=(0,\ell_e)$, whose image is the $1$-cell corresponding
to $e$, and the restriction $\map{\Phi_e}{I_e}{\Phi(I_e)} \subset X$
is a homeomorphism. The maps $\Phi_e$ induce a metric on $X$. In this
way, $X$ becomes a metric space.

Given a weighted discrete graph, we can abstractly construct the
associated metric graph as the disjoint union of the intervals $I_e$
for all $e \in E$ and together with a natural identification $\sim$ of
the end-points of these intervals (according to the combinatorial
structure of the graph), i.e.,
\begin{equation}
  \label{eq:constr.mg}
  X = \bigdcup_{e \in E} \clo I_e / {\sim}.
\end{equation}
We denote the union of the $0$-cells and the (disjoint) union of the
(open) $1$-cells (edges) by $X^0$ and $X^1$, respectively, i.e.,
\begin{equation*}
  X^0 = V \hookrightarrow X, \qquad 
  X^1 =\bigdcup_{e \in E} I_e \hookrightarrow X,
\end{equation*}
and both subspaces are canonically embedded in $X$.  

The metric graph $X$ becomes canonically a \emph{metric measure space}
by defining the distance of two points to be the length of the
shortest path in $X$, joining these points. We can think of the maps
$\map{\Phi_e}{I_e} X$ as coordinate maps and the Lebesgue measures on
the intervals $I_e$ induce a (Lebesgue) measure on the space $X$.

Since a metric graph is a topological space, and isometric to
intervals outside the vertices, we can introduce the notion of
measurability and differentiate function on the edges.  We start with
the basic Hilbert space
\begin{gather*}
  \Lsqr X := \bigoplus_{e \in E} \Lsqr {I_e}, \qquad
  f = \{f_e\}_e \quad \text{with $f_e \in \Lsqr{I_e}$ and}\\
  \normsqr f = \normsqr[\Lsqr X] f 
  := \sum_{e \in E} \int_{I_e} \abs{f_e(x)}^2 \dd x.
\end{gather*}

We define several types of Sobolev spaces on $X$. The \emph{maximal}
Sobolev space of order $k$ is given by
\begin{equation*}
  \Sobx[k] \max X
    := \bigoplus_{e \in E} \Sob [k]{I_e}
\end{equation*}
together with its natural norm.
The \emph{standard} or \emph{continuous} Sobolev space is given by
\begin{equation*}
  \Sob X := \Cont X \cap \Sobx \max X.
\end{equation*}
It can be shown that $\Sob X$ is indeed a Hilbert space as closed
subspace of the maximal Sobolev space using the length
condition~\eqref{eq:len.bd} (see e.g.~\cite[Lem.~5.2]{post:pre07a}).
For a graph with boundary $\bd V$, we define
\begin{equation*}
  \Sobx {\bd V} X := \bigset{f \in \Sob X}{f \restr {\bd V} = 0}
\end{equation*}
satisfying Dirichlet boundary conditions on $\bd V$. Again, $\Sobx
{\bd V} X$ is closed in $\Sobx \max X$. Note that $\Sobx V X =
\bigoplus_{e \in E} \Sobn {I_e}$ is the \emph{minimal} Sobolev space
of order $1$. We have the following inclusion of Sobolev spaces
\begin{equation}
  \label{eq:sob.incl}
  \Sobx V X 
  \subset \Sobx {\bd V} X
  \subset \Sob X
  \subset \Sobx \max X.
\end{equation}
We define quadratic forms in $\Lsqr X$ with domains
\begin{equation*}
  \dom \qf h^{\bd V} := \Sobx {\bd V} X, \qquad
  \dom \qf h  := \Sob X \und
  \dom \qf h^\Neu  := \Sobx \max X
\end{equation*}
acting as $\qf h^\bullet (f)= \normsqr{f'}=\sum_{e \in E} \int_{I_e}
\abs{f_e'}^2 \dd x$ in all cases. Denote by $\qlaplacian[\bd V]_X$,
$\qlaplacian_X$ and $\qlaplacian[\Neu]_X$ the corresponding
Laplacians, called \emph{Dirichlet(-Kirchhoff)}, \emph{Kirchhoff} and
\emph{fully decoupled Neumann} Laplacian. Note that
$\qlaplacian[\emptyset]_X = \qlaplacian_X$ and
\begin{equation*}
  \qlaplacian[V]_X = \bigoplus_{e \in E} \qlaplacian[\Dir]_{I_e} 
          \und 
  \qlaplacian[\Neu]_X = \bigoplus_{e \in E} \qlaplacian[\Neu]_{I_e}
\end{equation*}
are \emph{decoupled}, justifying the names \emph{fully decoupled
  Dirichlet} resp.\ \emph{Neumann Laplacian}.

A function $f$ is in the domain of the Dirichlet(-Kirchhoff) Laplacian
$\qlaplacian[\bd V]_X$ if and only if $f \in \Sobx[2] \max X$ and
\begin{subequations}
  \begin{gather}
    \label{eq:dir}
     f(v) = 0 \quad \forall \; v \in \bd V\\
    \label{eq:cont}
    \text{$f$ is continuous at each vertex $v \in \ring V=V \setminus
      \bd V$}\\
    \label{eq:kirch}
    \sum_{e \in E_v} f'_e(v)=0 \quad \forall \; v \in \ring V
  \end{gather}
\end{subequations}
where $f'_e(v)= -f_e'(0)$ if $v=\bd_-e$ and $f'_e(v)=f_e'(\ell_e)$
denotes the \emph{inwards} derivative of $f$ at the vertex $v$ along
the edge $e$.

If $X$ is a compact metric graph, the spectrum of all these operators
is purely discrete. We denote the eigenvalues by $\lambda_k^{\bd V}$,
$\lambda_k$ and $\lambda_k^\Neu$, $k=1,2,\dots$, respectively. It is
written in increasing order and respecting multiplicity.  Using the
variational characterisation of the eigenvalues, the \emph{min-max
  principle}~\eqref{eq:min.max} (see e.g.~\cite{davies:95}), we obtain
from the quadratic form inclusions~\eqref{eq:sob.incl} the 
reverse inequality for the corresponding eigenvalues, namely
\begin{equation*}
  \lambda_k^V \ge \lambda_k^{\bd V} \ge \lambda_k \ge \lambda_k^\Neu.
\end{equation*}

For an equilateral metric graph we obtain:
\begin{lemma}
  \label{lem:ew.edge}
  Assume that the metric graph $X$ is compact and all lengths $\ell_e$
  are equal to $1$, then
  \begin{equation*}
    (n+1)^2 \pi^2 
    = \lambda_k^V 
    \ge \lambda_k^{\bd V} 
    \ge \lambda_k 
    \ge \lambda_k^\Neu 
    = n^2 \pi^2
  \end{equation*}
  for $k=1 + n \abs{E},\dots, (n+1)\abs{E}$, $n=0,1,\dots$ In
  particular, the eigenvalues of the Dirichlet resp.\ Kirchhoff
  Laplacian on $X$ group into sets of cardinality $\abs{E}$
  (respecting multiplicity) lying inside the intervals
  $K_n:=[n^2\pi^2,(n+1)^2 \pi^2]$.
\end{lemma}

%
\section{Spectral relation between discrete and equilateral metric
  graphs}
\label{sec:sp.rel}
%
In this section, we give a complete description of the spectrum of the
standard discrete Laplacian and the Kirchhoff Laplacian (and the
corresponding Dirichlet versions on the boundary).  Outside the fully
decoupled Dirichlet spectrum $\Sigma^\Dir := \set{n^2
  \pi^2}{n=1,2,\dots}$, the relation is well-known, and there exist
more general results relating different spectral components also in
the case of infinite graphs (see
e.g.~\cite{von-below:85,exner:97b,cattaneo:97,pankrashkin:06a,post:pre07c,
  bgp:08} and the references therein).

Throughout this section, $G$ will denote a finite weighted graph with
standard weight $m(v)=\deg v$ and $m_e=1$. Moreover, $X$ will be the
associated compact metric graph with lengths $\ell_e=1$. We will refer
to such metric graphs also as \emph{equilateral}. To avoid
unnecessary exceptional cases, we assume that the graph is connected.
Some results hold also for non-compact graphs, see \Rem{non-compact}.

Denote by $\qlaplacian[\bd V]_X$ the metric graph Laplacian with Dirichlet
boundary conditions on $\bd V$ and Kirchhoff conditions on $\ring V$.
Similarly, let $\dlaplacian[\bd V]_G$ be the discrete Dirichlet Laplacian
associated to the underlying discrete graph $(G,\deg)$ with standard weights.
We denote by
\begin{equation*}
  \discr{\mc N}^{\bd V}(\eta) := \ker (\dlaplacian[\bd V]_G - \eta)
        \Und
  \mc N^{\bd V} (\lambda) := \ker (\qlaplacian[\bd V]_X - \lambda)
\end{equation*}
the corresponding eigenspaces.
\begin{proposition}
  \label{prp:sp.rel}
  Assume that the metric graph $X$ is compact and equilateral and set
  $\mu(\lambda):= 1 - \cos \sqrt \lambda$. Suppose in addition that
  $\lambda \notin \Sigma^\Dir$, i.e., $\mu(\lambda) \notin \{0,2\}$.
  Then the map
  \begin{equation*}
    \map {\Phi_\lambda} 
         {\discr{\mc N}^{\bd V}(\mu(\lambda))}
         {\mc N^{\bd V} (\lambda)}, \qquad
         F \mapsto f=\Phi_\lambda F
  \end{equation*}
  is an isomorphism where
  \begin{equation*}
    f_e(x) = F(\bd_-e) \frac{\sin \sqrt \lambda(1-x)}
                            {\sin \sqrt \lambda} +
             F(\bd_+e) \frac{\sin \sqrt \lambda x}
                            {\sin \sqrt \lambda}, \qquad \lambda > 0,
  \end{equation*}
  and $f_e(x)=F(\bd_-e) (1-x) + F(\bd_+e)x$ for $\lambda=0$.  In
  particular,
  \begin{equation*}
    \lambda \in \spec{\qlaplacian [\bd V]_X}
       \qquad \text{if and only if} \qquad
    \mu(\lambda) \in \spec{\dlaplacian[\bd V]_G}
  \end{equation*}
  (preserving the multiplicities of the eigenvalues).
\end{proposition}
The proof is straightforward. Note that it is the Kirchhoff boundary condition
leading to the discrete Laplacian expression and vice versa. The continuity
condition and the eigenvalue equation on the metric graph are automatically
fulfilled by this Ansatz.
\begin{definition}
  \label{def:ef.vx}
  We refer to the eigenfunctions $f=\Phi_\lambda F$ on the metric
  graph as \emph{(non-trivial) vertex-based} eigenfunctions, since
  they are completely determined by their values on the vertices and
  interpolated on the edges according to the solution of the
  differential equation.
\end{definition}

\subsection{Spectral relation at the Dirichlet spectrum}
The aim of the present subsection is to give a complete analysis of
the spectrum of $\qlaplacian [\bd V]_X$ at the exceptional values
$\lambda_n=n^2\pi^2 \in \Sigma^\Dir$. The multiplicity of these
eigenvalues was already calculated in~\cite{von-below:85} by a direct
proof not using the homology groups introduced in the next section.

We will show in the next lemma that there are two types of
eigenfunctions: the first type, vanishing at each vertex, is related
with the (relative) homology of the graph; the second type does not
vanish at any vertex and is related to the spectral points $0$ and $2$
of the discrete graph.
\begin{lemma}
  \label{lem:two.cases}
  Assume that $X$ is a connected compact equilateral metric graph and
  that $f \in \mc N^{\bd V}(\lambda_n)$. Then
  \begin{enumerate}
  \item either $f(v)=0$ for all vertices $v \in V$,
  \item or $f(v) \ne 0$ for all vertices $v \in V$. This case can only
    occur if there are no Dirichlet boundary conditions, i.e., $\bd V
    = \emptyset$.
  \end{enumerate}
  In the first case we have
  \begin{equation*}
    f_e(x)= \frac {f_e'(0)} {n\pi} \sin (n \pi x),
  \end{equation*}
  and in the latter case, $f$ is constant in all vertices if $n$ is even, or
  $f(\bd_+ e) = -f(\bd_-e)$ if $n$ is odd and $G$ is
  bipartite.
\end{lemma}
\begin{proof}
  Since $-f_e''= \lambda_n f_e$ on each edge, we must have
  \begin{equation}
    \label{eq:ansatz}
    f_e(x)= \alpha_e \cos(n \pi x) + \eta_e \sin (n \pi x).
  \end{equation}
  In particular, we have at a vertex $v=\bd_-e$ that $f(v)=f_e(0)=\alpha_e$
  and $f(v_e)=f_e(1)=\alpha_e (-1)^n$ and similarly if $v = \bd_+ e$.

  If $f(v)=0$ for a vertex $v$ then $\alpha_e=0$ hence also $f(v_e)=0$.  By
  the connectedness of the graph the first claim follows.

  If $f(v) \ne 0$, then $\alpha_e \ne 0$ and therefore $f(v_e)=(-1)^n
  f(v)$.  If $n$ is even, the second claim follows. The existence of a
  non-trivial function with alternating sign ($n$ odd) is an
  eigenfunction of the standard discrete Laplacian with eigenvalue $2$
  and therefore equivalent to the fact that the graph is bipartite
  (see \Prp{bipartite}).
\end{proof}

The previous lemma motivates the following definition:
\begin{definition}
  \label{def:ef.top}
  For the exceptional value $\lambda_n:= n^2\pi^2 \in \Sigma^\Dir$ ($n
  \ge 1$) we denote by
  \begin{equation*}
    {\mc N}_0^{\bd V}(\lambda_n) 
      := \bigset {f \in {\mc N}^{\bd V}(\lambda_n)}
                 {f(v)=0 \quad \forall \;v \in V}
  \end{equation*}
  the space of eigenfunctions vanishing at \emph{all} vertices. We
  call these eigenfunctions \emph{topological} or
  \emph{edge-based}.
\end{definition}
The name ``topological'' will be justified in \Sec{hom}, where we
relate this space with certain first homology groups. Note that these
eigenfunction still satisfy the Kirchhoff condition in the
\emph{inner} vertices which will give the relation with homology (see
especially \Prp{hom.ef}).

Let us state the following simple observation for general
eigenfunctions associated to $\lambda_n$:
\begin{lemma}
  \label{lem:kirch.hom}
  Assume that $f$ is written in the general form~\eqref{eq:ansatz}.
  Then $f$ fulfills the Kirchhoff conditions in all inner vertices $v
  \in \ring V$ iff $\eta=\{\eta_e\}_e \in \ker d_0^*$ if $n$ is even
  resp.\ $\eta \in \ker \unor d_0^*$ if $n$ is odd.
\end{lemma}
\begin{proof}
  From the form of $f$ on each edge, it follows $f_e'(0)=n \pi \eta_e$
  and $f_e'(1)=(-1)^n n \pi \eta_e$, i.e., the inwards derivative is
  given by $f_e'(v)= n \pi \orient \eta_e(v)$ if $n$ is even and
  $f_e'(v)=-n\pi \eta_e$ if $n$ is odd (recall that $f_e'(v)$ denotes
  the inward derivative, see \Sec{mg}). Now the Kirchhoff condition at
  $v \in \ring V$ is equivalent to $d_0^* \eta(v)=0$ resp.\ $\unor
  d_0^* \eta(v)=0$, since
  \begin{equation*}
    \sum_{e \in E_v} f_e'(v) 
    = n \pi \sum_{e \in E_v} \orient \eta_e(v)
    \Und
    \sum_{e \in E_v} f_e'(v)
    = - n \pi \sum_{e \in E_v} \eta_e
  \end{equation*}
  if $n$ is even or $n$ is odd, respectively.
  \end{proof}

\begin{definition}
  \label{def:ef.vx0}
  Assume that the graph $X$ is connected. In case that $n$ is odd, we
  assume furthermore that the graph is bipartite with corresponding
  partition $V = A \dcup B$, and that the graph is oriented such that
  $E=E^+(A,B)$, i.e., all edges start in $A$ and end in $B$. If $n$ is
  even, we do not need such an assumption.

  We call the function $\phi_n=\{\phi_{n,e}\}_e$ defined on each edge
  as
  \begin{equation*}
    \phi_{n,e}(x) = \cos (n \pi x)
  \end{equation*}
  the eigenfunction \emph{corresponding to the constant eigenfunction} if $n$
  is even and \emph{corresponding to the bipartite eigenfunction} if $n$ is
  odd. In both cases, we refer to $\phi_n$ as the \emph{trivial vertex-based
    eigenfunction}.
\end{definition}
The above names have the following justification: $\phi_n$ obviously
fulfills the eigenvalue equation for $\lambda_n$. Moreover, it
fulfills the Kirchhoff condition, since $\phi_{n,e}'(v)=0$ for all $e
\in E_v$.  If $n$ is even, then $\phi_{n,e}(0)=\phi_{n,e}(1)=1$, i.e.,
$\phi$ restricted to the vertices is the \emph{discrete constant}
eigenfunction.  If $n$ is odd, then the above defined function $\phi$
is continuous at each vertex, namely, $\phi_{n,e}(v)$ is independent
of $e \in E_v$. Moreover, $\phi_{n,e}(v)=\phi_{n,e}(0)=1$ if $v \in A$
and $\phi_{n,e}(v)=\phi_{n,e}(1)=-1$ if $v \in B$. In particular,
$F(v):=\phi_n(v)$ is the discrete bipartite eigenfunction. Note that
$F$ can be properly defined only in the bipartite case. Again, the
eigenfunction $\phi$ arises from a discrete eigenfunction, and is
interpolated on the edges, justifying the name ``vertex-based''
(cf.~\Def{ef.vx}).

We can express \Lem{two.cases} in terms of spaces:
\begin{proposition}
  \label{prp:dec.spaces}
  Assume that $X$ is a connected compact equilateral metric graph.
  \begin{enumerate}
  \item
    \label{case.1}
    If $\bd V \ne \emptyset$, then ${\mc N}^{\bd V}(\lambda_n) = {\mc
      N}_0^{\bd V}(\lambda_n)$.
  \item
    \label{case.2}
    If $\bd V=\emptyset$, then
    \begin{equation*}
      {\mc N}^{\bd V}(\lambda_n) = 
      \begin{cases}
        {\mc N}_0(\lambda_n) & \text{$n$ odd and $G$ not bipartite,}\\
        {\mc N}_0(\lambda_n) \oplus \C \phi_n & 
             \text{otherwise,}
      \end{cases}
    \end{equation*}
    where $\phi_n$ is defined in the previous definition.
  \end{enumerate}
\end{proposition}
\begin{proof}
  If $\bd V \ne\emptyset$ or if $\bd V = \emptyset$, $n$ is odd and
  $G$ is not bipartite, then \Lem{two.cases} implies that ${\mc
    N}^{\bd V}(\lambda_n) = {\mc N}_0^{\bd V}(\lambda_n)$.  This
  covers case~\eqref{case.1} and the first part of~\eqref{case.2}. In
  any other case there is, in addition to the space ${\mc N}_0^{\bd
    V}(\lambda_n)$, a trivial vertex-based eigenfunction $\phi_n$.  By
  the explicit characterisation of the elements in ${\mc N}_0^{\bd
    V}(\lambda_n)$ (cf.~\Def{ef.vx0} and \Lem{two.cases}) it is
  immediate that $\varphi_n$ is orthogonal to any function in ${\mc
    N}_0^{\bd V}(\lambda_n)$. This shows the first part in
  case~\eqref{case.2} and the proof is concluded.
\end{proof}

%
\section{Homology on graphs}
\label{sec:hom}
%
In order to understand the topological content of the eigenspace ${\mc
  N}_0^{\bd V}(\lambda_n)$, we introduce the concept of (relative)
homology for both, the oriented exterior derivative $d$ as well as for
the unoriented version $\unor d$. The main reason why we need both is
the fact, that in the case of even $n$, the function $x \mapsto
\sin(n\pi x)$ on an edge is \emph{antisymmetric} (with respect to the
middle point of $(0,1)$ and therefore encodes the orientation of the
edge. For odd $n$, the function is \emph{symmetric}, and the
orientation of an edge is irrelevant.  This material, in particular
the computation of the corresponding Betti numbers, will be crucial
for the eigenvalue bracketing in the next section and the relation
between metric and discrete eigenvalues.

Let $X$ be the topological graph associated to the \emph{finite} graph
$G$, and set $X^0=V$, $X^1=X \setminus X^0$. Then $X^1$ contains $\abs
E$-many components homeomorphic to $(0,1)$ and labelled by $e\in E$. Let
$C_p(X)$ be the group of $p$-chains with complex coefficients, i.e., the
vector space of formal sums
\begin{equation*}
  C_0(X) = \sum_{v \in V} \C \cdot v \und
  C_1(X) = \sum_{e \in E} \C \cdot e.
\end{equation*}
For a subset $\bd V$ of $V=X^0$ we define the group of \emph{relative}
$p$-chains as
\begin{equation*}
  C_p(X,\bd V) := C_p(X) / C_p(\bd V).
\end{equation*}
Note that since $\bd V$ consists only of points, we have the natural
identifications
\begin{equation*}
  C_0(X,\bd V)= C_0(\ring V)=\sum_{v \in \ring V} \C \cdot v \und
  C_1(X,\bd V)= C_1(X).
\end{equation*}

\subsection{Oriented homology}
The (oriented) boundary map $\map \bd {C_1(X)}{C_0(X)}$ is defined as
$\bd e = \bd_+ e - \bd_- e$, i.e., the formal difference of the
terminal minus the initial vertex of $e$. (We use the same symbol as
in the definition of the discrete graph since no confusion is
possible.)  In particular, for $c=\sum_{e \in E} \eta_e \cdot e$ we
have
\begin{align*}
  \bd c = \sum_{e \in E} \eta_e \cdot (\bd_+ e - \bd_-e)
        &= \sum_{v \in V} \Bigl(\sum_{e \in E^+_v} \eta_e
                              -\sum_{e \in E^-_v} \eta_e \Bigr) \cdot v\\
        &= \sum_{v \in V} \Bigl(\sum_{e \in E_v} 
                   \orient\eta_e(v)\Bigr) \cdot v
        = \sum_{v \in V} m(v) (d^* \eta)(v) \cdot v
\end{align*}
using~\eqref{eq:reorder} (recall that we assumed that $m_e=1$).  The
definition of the corresponding boundary map $\rel \bd$ is naturally given by
the commutativity of the diagram
\begin{equation*}
  \begin{diagram}
   0 & \rTo & C_0(\bd V) & \rTo & C_0(X) & \rTo & C_0(X,\bd V) & \rTo & 0\\
     &      &  \uTo^0    &      & \uTo^\bd &    & \uTo^{\rel \bd}&     &  \\
   0 & \rTo & C_1(\bd V) & \rTo & C_1(X) &\rTo  & C_1(X,\bd V) & \rTo & 0.
  \end{diagram}
\end{equation*}
In particular, we have 
\begin{equation*}
  \rel \bd e =
     \begin{cases}
           \bd_+ e - \bd_- e & \text{if $\bd_\pm e \in \ring V$,}\\
           \bd_+ e           & \text{if $\bd_+e \in \ring V$, 
                                   $\bd_-e \in \bd V$,}\\
           -\bd_- e          & \text{if $\bd_-e \in \ring V$, 
                                   $\bd_+e \in \bd V$}\\
           0                 & \text{if $\bd_\pm e \in \bd V$.}
  \end{cases}
\end{equation*}
Note that one can check as above that
\begin{equation}
  \label{eq:rel.bd}
  \rel \bd c = \sum_{v \in \ring V} m(v) (d_0^* \eta)(v) \cdot v.
\end{equation}
The corresponding homologies resp.\ relative homologies are now
defined as
\begin{align*}
  H_0(X) &:= C_0(X)/\ran \bd, & H_0(X,\bd V) &:= C_0(X,\bd V)/\ran \rel \bd\\
  H_1(X) &:= \ker \bd,        & H_1(X,\bd V) &:= \ker \rel \bd.
\end{align*}

\subsection{Unoriented homology}
The unoriented boundary map $\map {\unor \bd} {C_1(X)}{C_0(X)}$ is
defined similarly as
$\unor \bd e = \bd_+ e + \bd_- e$, i.e., the formal \emph{sum} of the
terminal and initial vertex of $e$. As before, we see that
\begin{equation*}
  \unor \bd c
        = \sum_{v \in V} m(v) (\unor d^* \eta)(v) \cdot v.
\end{equation*}
The corresponding unoriented relative boundary map is given as before
but just replacing $-\bd_-e$ by $+\bd_-e$.  Similarly, we have
\begin{equation}
  \label{eq:unor.rel.bd}
  \rel {\unor \bd} c = \sum_{v \in \ring V} m(v) (\unor d_0^* \eta)(v) \cdot v.
\end{equation}
The corresponding homologies resp.\ relative homologies are now
defined as
\begin{align*}
  \unor H_0(X) &:= C_0(X)/\ran \unor \bd, 
         & \unor H_0(X,\bd V) &:= C_0(X,\bd V)/\ran \rel {\unor \bd}\\
  \unor H_1(X) &:= \ker \unor \bd,        
         & \unor H_1(X,\bd V) &:= \ker \rel {\unor \bd}.
\end{align*}

\subsection{Calculation of the Betti numbers}
Denote by $b_p = b_p(X) = \dim H_p(X)$ the (oriented) Betti-numbers, and
similarly, $b_p^{\bd V} = b_p(X,\bd V)= \dim H_p(X,\bd V)$ the corresponding
relative Betti-numbers. Moreover, the corresponding notation with a bar, e.g.,
$\unor b_p = \dim \unor H_p(X)$ refers to the unoriented homology. The result
for the oriented Betti-numbers is standard. We include a short proof for the
unoriented case.
\begin{lemma}
  \label{lem:betti}
  Assume that the topological graph $X$ is compact and connected, and
  that $\bd V \ne \emptyset$. Then the oriented Betti numbers are
  given as
  \begin{align*}
    b_0(X)&= 1, &
    b_0(X,\bd V)&=0,\\
    b_1(X)&= \abs E - \abs V + 1, &
    b_1(X,\bd V)&=\abs E - \abs V + \abs{\bd V}.
  \end{align*}
  The unoriented Betti numbers are
  \begin{align*}
    \unor b_0(X)&= \beta &
    \unor b_0(X,\bd V)&= 0\\
    \unor b_1(X)&= \abs E - \abs V + \beta &
    \unor b_1(X,\bd V)&= \abs E - \abs V + \abs{\bd V}
  \end{align*}
  where $\beta=1$ if $X$ is bipartite and $0$ otherwise.
\end{lemma}
\begin{proof}
  We give the proof only for the unoriented case. It is more convenient to use
  the corresponding cohomologies, defined via
  \begin{equation*}
    \unor H^0(X) := \ker \unor d, \qquad
    \unor H^1(X) := \ker \unor d^*
  \end{equation*}
  and using the natural Hilbert space structure of the $\lsqrspace$-spaces
  with the standard weights $m(v)=\deg v$ and $m_e=1$.  Similarly, the
  relative cohomologies are defined as kernels of $\ker \unor d_0$ and $\ker
  \unor d_0^*$.  From~\eqref{eq:unor.rel.bd}, it is easy to see that the
  $p$-th relative homology and cohomology spaces are isomorphic,
  and similarly for the other cases.

  Moreover, $F \in \ker \unor d$ is equivalent to $0=\unor \dlaplacian
  F$, and by ~\eqref{eq:unor.lapl}, we conclude that $\dlaplacian F=2F$
  for the ``oriented'' Laplacian $\dlaplacian=d^* d$. Since $2$ is an
  eigenvalue of $\dlaplacian$ iff the graph is bipartite
  (cf.~\Prp{bipartite}), it follows that $\unor b_0(X)=\beta$ (recall
  that the graph is connected). The Euler characteristic is the same
  for the oriented and unoriented homology (see
  e.g.~\cite{post:pre07a}). Therefore $\unor b_1(X)=\unor b_0(X) -
  \chi(X) = \abs E - \abs V + \beta$.

  The relative Betti number $\unor b_0(X,\bd V)$ is easily seen to
  vanish, since the graph is connected and the function (the bipartite
  eigenfunction $F \in \ker (\dlaplacian - 2)$) is determined by its
  value at a single vertex.  To compute $\unor b_1(X,\bd V)$ we have
  to analyse $\ker \unor d_0^*$, where $\unor d_0^*=\iota^*\circ \unor
  d^*$ is given in \Sec{unor.ext.der}.  Note that
  \begin{equation*}
    \ker \unor d_0^* 
      = \ker \unor d^* \oplus 
        \bigset {\eta \in (\ker \unor d_0^*)^\orth}
                {\unor d^*\eta \restr {\ring V}=0}.
  \end{equation*}
  To compute the dimension of the second term of the previous equation note
  that
  \begin{align*}
   \dim \bigset{\eta \in (\ker \unor d_0^*)^\orth} 
               {\unor d^*\eta \restr {\ring V}=0}
    &= \dim \bigset{F \in \lsqr V}
                   {F \in \ran \unor d^*= (\ker \unor d)^\orth, \;
                        \supp F \subset \bd V }\\
    &= \abs{\bd V} - \beta.
  \end{align*}
  Altogether we have
  \begin{equation*}
    \unor b_1(X,\bd V) = \unor b_1(X)+ \abs{\bd V} - \beta
    = \abs E - \abs V + \abs{\bd  V}
  \end{equation*}
  and the proof is concluded.
\end{proof}

\subsection{The topological eigenspaces}
We can now relate the eigenfunctions vanishing at all vertices with
the homology. Recall that $\lambda_n:= n^2\pi^2$.
\begin{proposition}
  \label{prp:hom.ef}
  For any 1-chain $c=\sum_{e \in E} \eta_e \cdot e$ define $f_c \in
  \Lsqr X$ by $f_{c,e}(x):=\eta_e \sin (n \pi x)$. Then the mappings
  \begin{align*}
    \map {\Psi_n&} {H_1(X,\bd V)}
                   {{\mc N}_0^{\bd V}(\lambda_n)},
              & \text{$n \ne 0$ even, and}\\
    \map {\unor \Psi_n &} {\unor H_1(X,\bd V)}
                          {{\mc N}_0^{\bd V}(\lambda_n)},
              & \text{$n$ odd},
   \end{align*}
   given by $\Psi_n(c):=f_c$ and $\unor \Psi_n(c):=f_c$, respectively,
   are isomorphisms.
\end{proposition}
\begin{proof}
  We show first that $f_c \in {\mc N}_0^{\bd V}(\lambda_n)$. Note
  that, by construction $f_c\restr V=0$ and that $f_c$ is continuous
  on each vertex.  It remains to check the Kirchhoff condition at the
  inner vertices $v\in\ring V$.  Since $c=\sum_{e \in E} \eta_e \cdot
  e\in H_1(X,\bd V)$ we have that $\rel \bd c =0$, hence $d_0^*
  \eta=0$ with $\eta=(\eta_e)_e$ and $n\not=0$ even.
  From~\Lem{kirch.hom} we have that $f_c$ satisfies the Kirchhoff
  condition at $\ring V$.  Finally we have to show that $\Psi_n$ is
  bijective.  The injectivity of $\Psi_n$ is clear. In order to show
  the surjectivity, let $f \in {\mc N}_0^{\bd V}(\lambda_n)$ and put
  $\eta_e:= f'(0)/(n\pi)$. Then $\Psi_n(c)=f$ by construction, and
  $d_0^*\eta=0$. The case $n$ odd is done similarly.
\end{proof}
Note that for the topological eigenfunctions (or, what is the same,
edge-based) it is again the Kirchhoff condition giving the relation
with the discrete graph (or at least with its homology), as we have
already noticed for the vertex-based eigenfunctions in \Prp{sp.rel}.

\begin{remark}
  \label{rem:cattaneo}
  Note that Cattaneo~\cite{cattaneo:97} already calculated the
  spectrum of an equilateral (possibly infinite) graph (with $\bd
  V=\emptyset$) also for the exceptional values $\Sigma^\Dir$ without
  taking care about the multiplicities. She obtains the same result.
  Namely, if the graph has at least one even cycle (i.e., a closed
  path passing an even number of edges), then the first homology is
  non-trivial in the oriented and unoriented case ($b_1(X) \ge \unor
  b_1(X) >0$), and $\lambda_n$ is in the spectrum of $\qlaplacian_X$

  If $n$ is odd and the graph has only one odd cycle, then Cattaneo
  uses the following characterisation: $\lambda_n \in \qlaplacian_X$
  iff the graph is transient. The transience is equivalent to the
  existence of a \emph{flow} with finite energy and source $a$; in our
  notation, that there exists an element $\eta \in \lsqr E$ such that
  $d^* \eta = \delta_a$ ($\delta_a(v)=1$ if $a=v$ and $\delta_a(v)=0$
  otherwise).  The latter condition means that $\delta_a$ is in $\ran
  d^*$, i.e, orthogonal to $\ker d = \C\1_V$ if the graph is finite.
  But $\delta_a$ is never orthogonal to $\1_V$, so in this case, there
  are no eigenvalues, as we already conclude from $\unor b_1(X)=0$ and
  \Prp{hom.ef}.

  Note that Cattaneo's primary interest are Laplacians on infinite
  metric graphs with weights defined in a slightly different way than
  our metric graph Laplacians, see~\cite{cattaneo:97}.

  Moreover, von Below~\cite{von-below:85} already calculated the
  multiplicities of the exceptional eigenvalues $\lambda_n$ in the
  case $\bd V=\emptyset$, but without using homology groups.
\end{remark}

Although non-compact graphs are not our main purpose here, let us make
a few comment on this case. The non-compact case occurs in
\Secs{cov.gr}{ex} were we consider infinite covering graphs.
\begin{remark}
  \label{rem:non-compact}
  If $X$ is non-compact and connected, the spectral relation of
  \Prp{sp.rel} is still true, even more, one can show that all
  spectral types (discrete and essential, absolutely and singular
  continuous, (pure) point) are preserved, see~\cite{bgp:08} for
  details. Moreover, $\mc N^{\bd V}(\lambda_n)=\mc N_0^{\bd
    V}(\lambda_n)$ ($n \ge 1$) due to the fact that the trivial vertex
  based eigenfunctions $\phi_n$ are no longer in $\Lsqr X$. Moreover,
  we can easily extend the above results to the infinite case. Namely,
  if $n$ is even, \Prp{hom.ef} extends to the assertion that
  \begin{equation*}
    \map{\Psi^n}{H^1(X,\bd V)} {\mc N^{\bd V}(\lambda_n)},
    \qquad
    \eta \mapsto \sqrt 2 f_\eta, \qquad f_{\eta,e}(x)=\eta_e \sin(n\pi x)
  \end{equation*}
  is an isometric isomorphism using the corresponding
  $\lsqrspace$-cohomology $H^1(X,\bd)=\ker d_0^* \subset \lsqr E$.
  The case $n$ odd can be treated similarly.

\end{remark}

%
\section{Eigenvalue bracketing}
\label{sec:ev.brack}
%

\subsection{Eigenvalue counting for metric graphs}
\begin{figure}[h]
  \centering
\begin{picture}(0,0)%
  \includegraphics{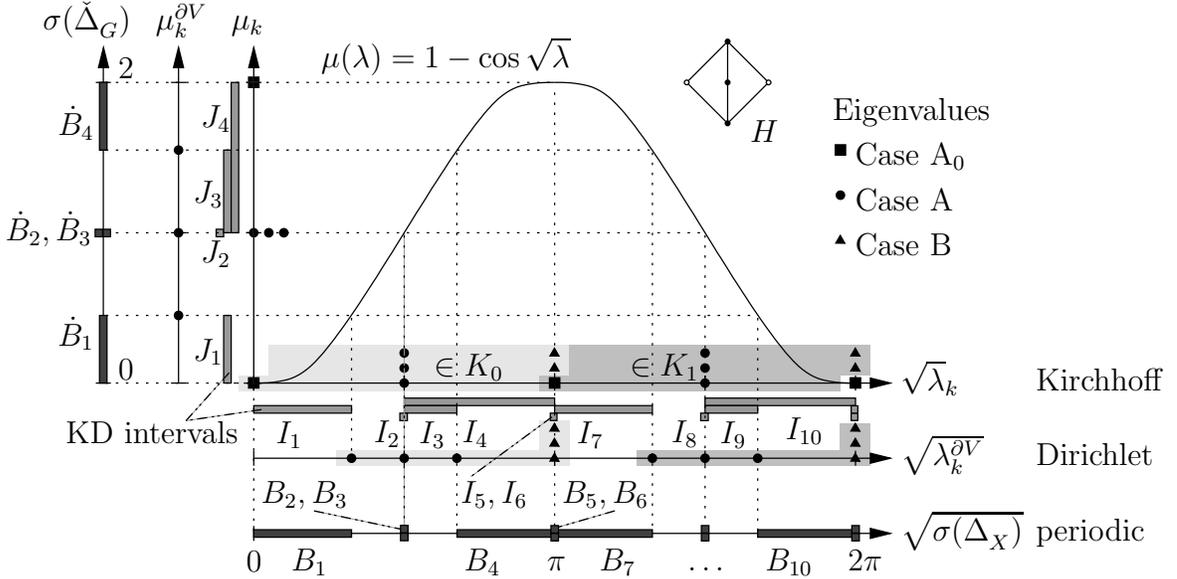}%
\end{picture}%
\setlength{\unitlength}{4144sp}%
\begin{picture}(6888,3421)(-1034,-3035)
\put(1531,-1816){$\in K_0$}
\put( 91,-1726){$J_1$}
\put(136,-1141){$J_2$}
\put( 91,-781){$J_3$}
\put(136,-331){$J_4$}
\put(2386,-2221){$I_7$}
\put(3634,-2198){$I_{10}$}
\put(3241,-2221){$I_9$}
\put(586,-2221){$I_1$}%
\put(1441,-2221){$I_3$}%
\put(1170,-2219){$I_2$}%
\put(2951,-2215){$I_8$}%
\put(1711,-2986){$B_4$}%
\put(3511,-2986){$B_{10}$}%
\put(2521,-2986){$B_7$}%
\put(1691,-2581){$I_5, I_6$}%
\put(1697,-2219){$I_4$}%
\put(2701,-1816){$\in K_1$}%
\put(4051,-556){Case A$_0$}%
\put(4051,-826){Case A}%
\put(4051,-1096){Case B}%
\put(4321,-1906){$\sqrt \lambda_k$}%
\put(-134,254){$\mu_k^{\bd V}$}%
\put(4321,-2356){$\sqrt {\lambda_k^{\bd V}}$}%
\put(-809,254){$\spec{\dlaplacian_G}$}%
\put(-359,-39){$2$}%
\put(-359,-1840){$0$}%
\put(4321,-2806){$\sqrt{\spec{\qlaplacian_X}}$}%
\put(2206,-2986){$\pi$}%
\put(406,-2991){$0$}%
\put(4006,-2986){$2\pi$}%
\put(3916,-286){Eigenvalues}%
\put(856, 29){$\mu(\lambda)=1- \cos \sqrt \lambda$}%
\put(5131,-1906){Kirchhoff}%
\put(5131,-2356){Dirichlet}%
\put(5131,-2806){periodic}%
\put(-671,-2221){KD intervals}%
\put(3421,-421){$H$}%
\put(316,254){$\mu_k$}%
\put(676,-2986){$B_1$}%
\put(496,-2581){$B_2, B_3$}%
\put(-719,-1636){$\dot B_1$}%
\put(-719,-376){$\dot B_4$}%
\put(-1034,-1006){$\dot B_2, \dot B_3$}%
\put(2296,-2581){$B_5, B_6$}%
\put(3049,-2986){\dots}%
\end{picture}%
  \caption{The various eigenvalues for the bipartite graph with
    fundamental domain $H$ and periodic graph $G$ of \Fig{gr.bip} with
    five vertices, two boundary vertices and six edges. Multiple
    eigenvalues are indicated by repeated symbols, compare with the
    tables in \Sec{ev.brack}. The eigenvalues for the Kirchhoff and
    Dirichlet metric Laplacian are grouped into members of six (light
    grey and dark grey) belonging to $K_0=[0,\pi^2]$ and
    $K_1=[\pi^2,4\pi^2]$ as predicted in \Lem{ew.edge}.  For a
    discussion of the relation of the KD intervals with periodic
    operators see \Ex{bip}.}
  \label{fig:sp.rel}
\end{figure}

Let us now combine the results of the previous sections. In
particular, we will show how the $\abs V$ eigenvalues $\mu_k$ of the
discrete Laplacian are related with the $\abs E$ eigenvalues
$\lambda_k$ in $K_n=[n^2\pi^2, (n+1)^2\pi^2]$ of the Kirchhoff
Laplacian. For the Dirichlet operators we relate the $\abs V -
\abs{\bd V}$ discrete eigenvalues $\mu_k^{\bd V}$ with $\abs E$ metric
eigenvalues $\lambda_k^{\bd V} \in K_n$.  In
\Figs{sp.rel}{sp.rel.n-bp} we illustrated the spectral relations for a
bipartite and non-bipartite graph of \Exs{bip}{non-bip} (see
\Figs{gr.bip}{gr.non-bip}).  Doing a neat bookkeeping one can check
the different possibilities given in the tables below.

We start with a basic definition:
\begin{definition}
  \label{def:kd.int.mg}
  We define the (metric) Kirchhoff-Dirichlet intervals $I_k=I_k(X,\bd V)$ of
  the metric graph $X$ with boundary $\bd V$ as
  \begin{equation*}
    I_k := [\lambda_k, \lambda_k^{\bd V}], \qquad k=1,2,\dots
  \end{equation*}
\end{definition}
Note that by \Lem{ew.edge}, the interval is non-empty and $I_k \subset
K_n$ for $k=n \abs E +1, \dots, (n+1) \abs E$, where
$K_n:=[n^2\pi^2,(n+1)^2\pi^2]$ for $n=0,1,\dots$

The aim of the following eigenvalue counting is to understand the
nature of the intervals $I_k$, i.e., whether they reduce to points or
are contained in $\ring K_n$. It is therefore unavoidable to give a
precise account of the eigenvalues repeated according to multiplicity
in the Kirchhoff as well as in the Dirichlet case, distinguishing
bipartite and non-bipartite graphs.
 
\begin{figure}[h]
  \centering
\begin{picture}(0,0)%
  \includegraphics{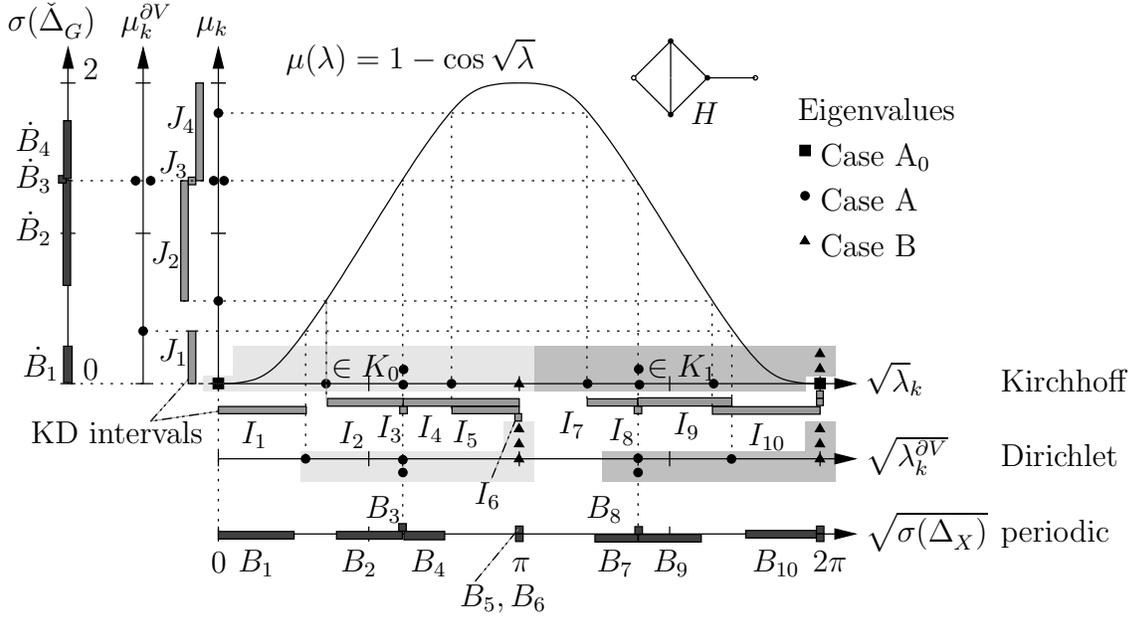}%
\end{picture}%
\setlength{\unitlength}{4144sp}%
\begin{picture}(6663,3626)(-809,-3240)
\put(1126,-1816){$\in K_0$}%

\put(586,-2221){$I_1$}%

\put(1170,-2219){$I_2$}%

\put( 91,-1726){$J_1$}%

\put( 46,-1141){$J_2$}%

\put( 91,-601){$J_3$}%

\put(136,-331){$J_4$}%

\put(1396,-2181){$I_3$}%

\put(1981,-2581){$I_6$}%

\put(1639,-2184){$I_4$}%

\put(3172,-2166){$I_9$}%

\put(2785,-2201){$I_8$}%

\put(2486,-2159){$I_7$}%

\put(3620,-2230){$I_{10}$}%

\put(1851,-2214){$I_5$}%

\put(3016,-1816){$\in K_1$}%

\put(4051,-556){Case A$_0$}%

\put(4051,-826){Case A}%

\put(4051,-1096){Case B}%

\put(4321,-1906){$\sqrt \lambda_k$}%

\put(-134,254){$\mu_k^{\bd V}$}%

\put(4321,-2356){$\sqrt {\lambda_k^{\bd V}}$}%

\put(-809,254){$\spec{\dlaplacian_G}$}%

\put(-359,-39){$2$}%

\put(-359,-1840){$0$}%

\put(4321,-2806){$\sqrt{\spec{\qlaplacian_X}}$}%

\put(2206,-2986){$\pi$}%

\put(406,-2991){$0$}%

\put(4006,-2986){$2\pi$}%

\put(3916,-286){Eigenvalues}%

\put(856, 29){$\mu(\lambda)=1- \cos \sqrt \lambda$}%

\put(5131,-1906){Kirchhoff}%

\put(5131,-2356){Dirichlet}%

\put(5131,-2806){periodic}%

\put(-671,-2221){KD intervals}%

\put(3277,-306){$H$}%

\put(316,254){$\mu_k$}%

\put(-764,-1006){$\dot B_2$}%

\put(-719,-1816){$\dot B_1$}%

\put(-764,-691){$\dot B_3$}%

\put(-764,-466){$\dot B_4$}%

\put(570,-2986){$B_1$}%

\put(1184,-2986){$B_2$}%

\put(1330,-2667){$B_3$}%

\put(1896,-3195){$B_5, B_6$}%

\put(2714,-2986){$B_7$}%

\put(3644,-2986){$B_{10}$}%

\put(3049,-2986){$B_9$}%

\put(1598,-2986){$B_4$}%

\put(2650,-2667){$B_8$}%

\end{picture}%
  \caption{The various eigenvalues for the non-bipartite graph with
    fundamental domain $H$ and periodic graph $G$ of \Fig{gr.non-bip} with
    five vertices, two boundary vertices and six edges. Again,
    multiple eigenvalues are indicated by repeated symbols; and the
    eigenvalues for the Kirchhoff and Dirichlet metric Laplacian are
    grouped into members of six (light grey and dark grey) as in
    \Fig{sp.rel} For a discussion of the relation of the KD intervals
    with periodic operators see \Ex{non-bip}.}
  \label{fig:sp.rel.n-bp}
\end{figure}
\subsubsection*{Counting the Kirchhoff eigenvalues}
The following result summarises several facts of this and the previous section.
In particular it is a consequence of
\Prpss{sp.rel}{dec.spaces}{hom.ef}. We use the abbreviations EF for
eigenfunction and EV for eigenvalue. The trivial vertex-based
eigenfunction $\phi_n$ is described in \Def{ef.vx0}, the non-trivial
vertex-based eigenfunctions are described in \Prp{sp.rel} (see
\Def{ef.vx}) and the topological eigenfunctions are described in
\Prp{hom.ef} (see \Def{ef.top}).
\begin{proposition}
  \label{prp:count.kirch}
  Let $X$ be a connected compact equilateral metric graph and let
  $n=0,1,2,\dots$ The EVs $\lambda_k$ of $\qlaplacian_X$ distribute in
  groups of $\abs E$ EVs contained in the intervals $K_n:=[n^2\pi^2,
  (n+1)^2\pi^2]$. We list them in the following tables according to
  the various possibilities.  The brace under the range of the index
  $k$ denotes the number of such indices.

  If the graph is bipartite we have:
  \begin{center}
    \begin{tabular}{|l|c|l|p{12ex}|l|}
      \multicolumn{5}{c}{\textbf{Metric Kirchhoff eigenvalues 
                  for bipartite graphs}}\\
      \hline
      Case & Range of index $k$ & $\lambda_k$ & Type of EF &
      EF described in\\
      \hline\hline
      A$_0$ & $\underbrace{n \abs E + 1}_{1}$ & 
      $= n^2\pi^2$ & $\phi_n$, trivial\newline vertex-based& 
           \Prp{dec.spaces}\\ \hline
      A & $\underbrace{n \abs E + 2, \dots, n\abs E + \abs V - 1}_{\abs V -
        2}$  & $\in \ring K_n$
      & vertex-based& \Prp{sp.rel}\\ \hline
      B & $\underbrace{n\abs E + \abs V, \dots, (n+1) \abs
        E}_{b_1(X)=
         \unor b_1(X)=\abs E - \abs V + 1}$  
      & $=(n+1)^2 \pi^2$
      & topological& \Prp{hom.ef}\\ \hline
    \end{tabular}
  \end{center}

  If the graph is not bipartite we have:
  \begin{center}
    \begin{tabular}{|l|c|l|p{12ex}|l|}
      \multicolumn{5}{c}{\textbf{Metric Kirchhoff eigenvalues for
          non-bipartite graphs, $n$ even}}\\
      \hline
      Case & Range of index $k$ & $\lambda_k$ & Type of EF &
      EF described in\\
      \hline\hline
      A$_0$ & $\underbrace{n \abs E + 1}_{1}$ & 
      $= n^2\pi^2$ & $\phi_n$, trivial \newline 
           vertex-based& \Prp{dec.spaces}\\ \hline
      A & $\underbrace{n \abs E + 2, \dots, n\abs E + \abs V}_{\abs V -
        1}$  & $\in \ring K_n$
      &  vertex-based& \Prp{sp.rel}\\ \hline
      B & $\underbrace{n\abs E + \abs V+1, \dots, (n+1) \abs
        E}_{\unor b_1(X)=\abs E - \abs V}$  & $=(n+1)^2 \pi^2$
      & topological& \Prp{hom.ef}\\ \hline
    \end{tabular}
  \end{center}

  \begin{center}
    \begin{tabular}{|l|c|l|p{12ex}|l|}
      \multicolumn{5}{c}{\textbf{Metric Kirchhoff eigenvalues for
          non-bipartite graphs, $n$ odd}}\\
      \hline
      Case & Range of index $k$ & $\lambda_k$ & Type of EF &
      EF described in\\
      \hline\hline
      A & $\underbrace{n \abs E + 1, \dots, n\abs E + \abs V-1}_{\abs V -
        1}$  & $\in \ring K_n$
      & non-trivial \newline vertex-based& \Prp{sp.rel}\\ \hline
      B & $\underbrace{n\abs E + \abs V, \dots, (n+1) \abs
        E}_{b_1(X)=\abs E - \abs V +1}$  & $=(n+1)^2 \pi^2$
      & topological& \Prp{hom.ef}\\ \hline
    \end{tabular}
  \end{center}
\end{proposition}
\begin{remark}
  \samepage  
  \indent
  \begin{enumerate}
  \item For a bipartite graph, the trivial vertex-based eigenfunction
    $\phi_n$ corresponds to the constant discrete EF if $n$ is even
    and to the bipartite eigenfunction if $n$ is odd.
  \item Note that in the non-bipartite case, there is one eigenvalue
    of Case~A more than in the bipartite case. In the bipartite case,
    this additional eigenfunction is either a topological one (Case~B)
    if $n$ is even or a trivial vertex-based one (Case~A) if $n$ is
    odd, namely the one corresponding to the bipartite EF.
  
  \end{enumerate}
\end{remark}
\subsubsection*{Counting the Dirichlet eigenvalues}
The Dirichlet case is simpler and does not distinguish
the bipartite and non-bipartite case.
\begin{proposition}
  \label{prp:count.dir}
  Let $X$ be a connected compact equilateral metric graph with non-trivial
  boundary $\bd V \ne \emptyset$ and let $n=0,1,2,\dots$ The EVs $\lambda^{\bd
    V}_k$ of $\qlaplacian[\bd V] X$ distribute in groups of $\abs E$ EVs
  contained in the intervals $K_n:=[n^2\pi^2, (n+1)^2\pi^2]$. We list them in
  the following table:
  \begin{center}
    \begin{tabular}{|l|c|l|p{12ex}|}
      \multicolumn{4}{c}{\textbf{Metric Dirichlet eigenvalues}}\\
      \hline
      Case & Range of index $k$ & $\lambda^{\bd V}_k$ & 
      Type of EF\\
      \hline\hline
      A & $\underbrace{n \abs E + 1, \dots, n\abs E + \abs V -
        \abs{\bd V}}_{\abs V - \abs{\bd V}}$  & $\in \ring K_n$
      & non-trivial \newline vertex-based\\ \hline
      B & $\underbrace{n\abs E + \abs V - \abs{\bd V} + 1, \dots, (n+1) \abs
        E}_{b_1(X,\bd V)=
        \unor b_1(X,\bd V)=\abs E - \abs
        V + \abs{\bd V}}$  
      & $=(n+1)^2 \pi^2$
      & topological\\ \hline
    \end{tabular}
  \end{center}
  Again, Case~A is described in \Prp{sp.rel} and Case~B in \Prp{hom.ef}.
\end{proposition}

We can now describe precisely all possible combinations of
Kirchhoff-Dirichlet intervals that arise from the previous tables:
\begin{proposition}
  \label{prp:kd-int.mg}
  Let $X$ be a connected compact equilateral metric graph with
  non-empty boundary $\bd V\ne \emptyset$ and let $n=0,1,2,\dots$ The
  metric Kirchhoff-Dirichlet intervals are given in the table below.
  We call an interval \emph{non-degenerate} if its interior is
  non-empty. The case-labeling refers to the cases of the Kirchhoff
  (first letter) and Dirichlet (second letter) eigenvalue.

  If the graph is bipartite we have:
  \begin{center}
    \begin{tabular}{|l|c|l|p{17ex}|}
      \multicolumn{4}{c}
      {\textbf{Metric Kirchhoff-Dirichlet intervals for bipartite graphs}}\\
      \hline
      Case & Range of index $k$ & $I_k=I_k(X,\bd V)$ & 
      Type of interval\\
      \hline\hline
      A$_0$A & $\underbrace{n \abs E + 1}_{1}$ & 
      $= [n^2\pi^2,\lambda_k^{\bd V}]$ & non-degenerate\\
      \hline
      AA & $\underbrace{n \abs E + 2, \dots, n\abs E + \abs V 
        - \abs{\bd V}}_{\abs V - \abs{\bd V} - 1}$ 
      & $\subset \ring K_n$
      & degenerate or \newline non-degenerate\\ \hline
      AB & $\underbrace{n \abs E + \abs V -  \abs{\bd V} + 1, 
        \dots, n\abs E + \abs V - 1}_{\abs{\bd V} - 1}$ 
      & $=[\lambda_k, (n+1)^2\pi^2]$
      & non-degenerate\\ \hline
      BB & $\underbrace{n\abs E + \abs V, \dots, (n+1) \abs
        E}_{b_1(X) =
        \unor b_1(X)=\abs E - \abs V + 1}$  
      & $=(n+1)^2 \pi^2$
      & degenerate\\ \hline
    \end{tabular}
  \end{center}

  For non-bipartite graphs, we obtain:
  \begin{center}
    \begin{tabular}{|l|c|l|p{17ex}|}
      \multicolumn{4}{c}
      {\textbf{Metric Kirchhoff-Dirichlet intervals for 
                      non-bipartite graphs ($n$ even)}}\\
      \hline
      Case & Range of index $k$ & $I_k=I_k(X,\bd V)$ & 
      Type of interval\\
      \hline\hline
      A$_0$A & $\underbrace{n \abs E + 1}_{1}$ & 
      $= [n^2\pi^2,\lambda_k^{\bd V}]$ & non-degenerate\\
      \hline
      AA & $\underbrace{n \abs E + 2, \dots, n\abs E + \abs V 
        - \abs{\bd V}+1}_{\abs V - \abs{\bd V}}$ 
      & $\subset \ring K_n$
      & degenerate or \newline non-degenerate\\ \hline
      AB & $\underbrace{n \abs E + \abs V -  \abs{\bd V} + 2, 
        \dots, n\abs E + \abs V}_{\abs{\bd V}}$ 
      & $=[\lambda_k, (n+1)^2\pi^2]$
      & non-degenerate\\ \hline
      BB & $\underbrace{n\abs E + \abs V + 1, \dots, (n+1) \abs
        E}_{\unor b_1(X)=\abs E - \abs V}$  
      & $=(n+1)^2 \pi^2$
      & degenerate\\ \hline
    \end{tabular}
  \end{center}
  \begin{center}
    \begin{tabular}{|l|c|l|p{17ex}|}
      \multicolumn{4}{c}
      {\textbf{Metric Kirchhoff-Dirichlet intervals for 
                   non-bipartite graphs ($n$ odd)}}\\
      \hline
      Case & Range of index $k$ & $I_k=I_k(X,\bd V)$ & 
      Type of interval\\
      \hline\hline
      AA & $\underbrace{n \abs E + 1, \dots, n\abs E + \abs V 
        - \abs{\bd V}}_{\abs V - \abs{\bd V}}$ 
      & $\subset \ring K_n$
      & degenerate or \newline non-degenerate\\ \hline
      AB & $\underbrace{n \abs E + \abs V -  \abs{\bd V} + 1, 
        \dots, n\abs E + \abs V - 1}_{\abs{\bd V} - 1}$ 
      & $=[\lambda_k, (n+1)^2\pi^2]$
      & non-degenerate\\ \hline
      BB & $\underbrace{n\abs E + \abs V, \dots, (n+1) \abs
        E}_{b_1(X) =
        \unor b_1(X)=\abs E - \abs V + 1}$  
      & $=(n+1)^2 \pi^2$
      & degenerate\\ \hline
    \end{tabular}
  \end{center}
\end{proposition}

\subsection{Eigenvalue counting for discrete graphs}
We can now carry over the eigenvalue monotonicity of the Kirchhoff and
Dirichlet \emph{metric} Laplacian to the \emph{discrete} one.  For
this purpose it is is enough to consider only the metric graph
eigenvalues in the \emph{first} interval $K_0=[0,\pi^2]$, since on
this interval, the function $\mu(\lambda)=1-\cos(\sqrt \lambda)$ is
increasing. Denote by $\mu_k$ ($k=1,\dots, \abs V$) the eigenvalues of
the standard discrete Laplacian $\dlaplacian_G$, and by $\mu_k^{\bd
  V}$ ($k=1,\dots, \abs V -\abs {\bd V}$) the eigenvalues of the
(standard) discrete Laplacian with Dirichlet conditions on $\bd V$
(see \Sec{dg}), in both cases counted according to multiplicity.
\begin{definition}
  \label{def:kd.int.dg}
  We define the (discrete) Kirchhoff-Dirichlet intervals $J_k=J_k(G,\bd
  V)$ of the metric graph $X$ with boundary $\bd V$ as
  \begin{equation*}
    J_k := [\mu_k, \mu_k^{\bd V}], \qquad k=1,2,\dots, \abs V - \abs{\bd V}.
  \end{equation*}
  For higher indices, we set
  \begin{equation*}
    J_k := [\mu_k, 2], \qquad 
          k = \abs V - \abs{\bd V} + 1, \dots, \abs V.
  \end{equation*}
\end{definition}
\begin{remark}
  Note that the names ``Kirchhoff'' and ``Dirichlet'' for the standard
  Laplacian (with Dirichlet conditions on $\bd V$) is justified by
  \Prp{sp.rel}. Note also, that the operators act in spaces of
  different dimensions. In particular, the standard Laplacian with
  $\bd V=\emptyset$ can be written as a $\abs V \times \abs V$-matrix
  and has therefore $\abs V$ eigenvalues. Similarly, the standard
  Dirichlet Laplacian has $\abs V - \abs{\bd V}$ eigenvalues.
\end{remark}
{}From \Prp{kd-int.mg} we immediately obtain:
\begin{proposition}
  \label{prp:kd-int.dg}
  Let $G$ be a connected finite discrete graph with standard weight
  ($m(v)=\deg v$ and $m_e=1$) and non-trivial boundary $\bd V\ne
  \emptyset$. Then the discrete Kirchhoff-Dirichlet intervals are
  given in the table below. (Note that the type of the interval is
  the same as the type for the metric graph.)

  If the graph is bipartite we have:
  \begin{center}
    \begin{tabular}{|l|c|l|p{17ex}|}
      \multicolumn{4}{c}
      {\textbf{Discrete Kirchhoff-Dirichlet intervals for bipartite graphs}}\\
      \hline
      Case & Range of index $k$ & $J_k=J_k(G,\bd V)$ & 
      Type of interval\\
      \hline\hline
      A$_0$A & $k=1$ & 
      $= [0, \mu_k^{\bd V}]$ & non-degenerate\\
      \hline
      AA & $\underbrace{2, \dots, \abs V 
        - \abs{\bd V}}_{\abs V - \abs{\bd V} - 1}$ 
      & $\subset (0,2)$
      & degenerate or \newline non-degenerate\\ \hline
      AB & $\underbrace{\abs V -  \abs{\bd V} + 1, 
        \dots, \abs V - 1}_{\abs{\bd V} - 1}$ 
      & $=[\mu_k, 2]$
      & non-degenerate\\ \hline
      BB & $k=\abs V$ 
      & $=\{2\}$
      & degenerate\\ \hline
    \end{tabular}
  \end{center}
  For non-bipartite graphs, we obtain:
  \begin{center}
    \begin{tabular}{|l|c|l|p{17ex}|}
      \multicolumn{4}{c}
      {\textbf{Discrete Kirchhoff-Dirichlet intervals for 
          non-bipartite graphs}}\\
      \hline
      Case & Range of index $k$ & $J_k=J_k(G,\bd V)$ & 
      Type of interval\\
      \hline\hline
      A$_0$A & $k=1$ & 
      $= [0, \mu_k^{\bd V}]$ & non-degenerate\\
      \hline
      AA & $\underbrace{2, \dots, \abs V 
        - \abs{\bd V}+1}_{\abs V - \abs{\bd V}}$ 
      & $\subset (0,2)$ & degenerate or \newline non-degenerate\\ \hline
      AB & $\underbrace{\abs V -  \abs{\bd V} + 2, 
        \dots, \abs V}_{\abs{\bd V}}$ 
      & $=[\mu_k, 2]$
      & non-degenerate\\ \hline
    \end{tabular}
  \end{center}
\end{proposition}

\subsection{Spectral symmetry for bipartite graphs}
Let us carry over the spectral symmetry for \emph{discrete} bipartite
graphs already mentioned in \Prp{bipartite} to the metric case. Note
that the symmetry function in the discrete case is
\begin{equation*}
  \map \theta {[0,2]}{[0,2]}, \qquad \theta(\mu)=2-\mu.
\end{equation*}
In particular, the fixed point of $\theta$, i.e., $\mu=1$, is
always an eigenvalue of $\dlaplacian[\bd V]_G$ if $\abs V - \abs {\bd
  V}$ is odd.

We recall the definition $K_n:=[n^2\pi^2,(n+1)^2\pi^2]$.
\begin{proposition}
  \label{prp:bip.mg}
  Suppose that $X$ is a bipartite equilateral compact metric graph
  with boundary $\bd V$ ($\bd V$ may be empty) and let $\lambda \in
  \ring K_n$. Then
  \begin{equation*}
    \lambda \in \spec {\qlaplacian[\bd V]_X } \qquad \text{iff} \qquad
    \tau_n(\lambda) \in \spec {\qlaplacian[\bd V]_X }
  \end{equation*}
  and the multiplicity is preserved. Here,
  \begin{equation*}
    \map{\tau_n}{K_n}{K_n}, \qquad
    \tau_n(\lambda) := \bigl((2n+1)\pi - \sqrt \lambda\bigr)^2.
  \end{equation*}
  If $\abs V - \abs {\bd V}$ is odd, then the fixed point of $\tau_n$,
  i.e., $\lambda=(n+1/2)^2 \pi^2$ is an eigenvalue of $\qlaplacian[\bd
  V]_X$.

  Moreover, for $n \ge 1$ the map $\tau_n$ also interchanges the
  topological eigenvalues $\lambda_n=n^2 \pi^2$ and
  $\lambda_{n+1}=(n+1)^2 \pi^2$. The corresponding eigenfunctions for
  $\tau_n(\lambda_n) = \lambda_{n+1}$ are obtained by those from
  $\lambda_n$ by keeping the amplitude of the oscillation and interchanging the
  frequency.

  If $X$ is non-compact, the spectral symmetry
  $\tau_n(\spec{\qlaplacian_X} \cap K_n) = \spec{\qlaplacian_X} \cap
  K_n$ still holds.
\end{proposition}
\begin{proof}
  The results for eigenvalues $\lambda \in \ring K_n$ follow immediately from
  \Prps{bipartite}{sp.rel}. Also, the trivial vertex-based eigenfunctions are
  interchanged by the symmetry, as in the discrete case. For the topological
  eigenvalues, note that their structure is given in \Prp{hom.ef}, and that
  the oriented and unoriented Betti numbers agree, namely $b_1(X,\bd V)=\unor
  b_1(X,\bd V)$. The non-compact case follows by the spectral relation for
  $\lambda \in \ring K_n$ (see \Rem{non-compact}), and by the closeness of the
  spectrum for the endpoints of $K_n$.
\end{proof}
%
\section{Equivariant Laplacians and coverings}
\label{sec:eq.lapl}
%

In the sequel, we will analyse metric and discrete Laplacians on
covering graphs.

\subsection{Equivariant metric Laplacians}
We start with a \emph{metric} covering graph $X \to X_0$ with covering
group $\Gamma$ (in general non-abelian) and compact quotient $X_0$,
see also \cite[Sec.~6]{sunada:pre08} for related aspects. We call the
metric Laplacian $\qlaplacian_X$ on $X$ also \emph{$\Gamma$-periodic}.
A \emph{fundamental domain} of a metric graph covering $X \to X_0$ is
a closed subset $Y$ of $X$ such that
\begin{equation*}
  \gamma \ring Y \cap \ring Y = \emptyset, \quad \gamma\ne 1, \qquad
  \bigcup_{\gamma \in \Gamma} \gamma Y = X.
\end{equation*}
Note that the interior of a fundamental domain $\ring Y$ can always be
embedded isometrically in the quotient graph $X_0$.  Moreover, we
assume that the boundary of $Y$ (as topological subset of $X$)
consists only of vertices, which are precisely the \emph{boundary}
vertices, i.e.,
\begin{equation}
  \label{eq:fd.bd}
  \bd V := \bd Y=Y\setminus\ring Y \subset V.
\end{equation}
Since we can interpret $\ring Y$ as subset of $X_0$, we define the set
of \emph{inner vertices} of the quotient $X_0$ by $\ring V_0 := \ring
Y \cap V_0$, depending of course on the fundamental domain.

Associated to a fundamental domain is a metric graph also denoted by
the symbol $Y$ with boundary vertices $\bd V=\bd Y$ (\emph{not}
embedded in the quotient). We define the Dirichlet and Kirchhoff
metric Laplacians on this graph, namely, we consider $\qlaplacian[\bd
V]_{Y}$ and $\qlaplacian_{Y}$ defined via their quadratic forms on
$\Sobx {\bd V} {Y}$ and $\Sob {Y}$.

Let $\rho$ be a unitary representation of $\Gamma$, i.e., $\rho$ is a
homomorphism from $\Gamma$ into the group of unitary operators on some
Hilbert space $\HS$.  In order to analyse the spectrum of the periodic
operator, we need the following definition:
\begin{definition}
  \label{def:eq.lapl}
  A function $\map f X \HS$ is called \emph{equivariant} iff
  \begin{equation*}
    f(\gamma \cdot x)= \rho(\gamma) f(x), \qquad \forall\; 
                x \in X, \gamma \in \Gamma.
  \end{equation*}
\end{definition}
Clearly, a $\rho$-equivariant function, locally in $\Sobspace 1$, is
determined by its values on $Y$, namely, the equivariance
condition \Def{eq.lapl} reduces to a condition for the boundary
vertices $x=v \in \bd V$ such that $\gamma \cdot v \in \bd V$. We therefore
set
\begin{equation*}
  \Sobx \rho {X_0,\HS} := 
     \bigset{f \in \Sob {Y} \otimes \HS}
       {f(\gamma \cdot v)=\rho(\gamma) f(v) \quad 
             \forall\; v \in \bd V \text{ such that } 
                     \gamma \cdot v \in \bd V}.
\end{equation*}
We can consider functions in $\Sobx \rho {X_0, \HS}$ as functions on
the quotient metric graph $X_0$, where the continuity condition at the
boundary vertices is replaced by the equivariance condition.

Denote by $\qlaplacian[\rho]_{X_0}$ the operator associated to the
quadratic form
\begin{equation}
  \label{eq:def.qf}
  \qf h^\rho(f):= 
  \sum_{e \in E_0} \int_0^{\ell_e} \normsqr[\HS] {f'(x)} \dd x,
  \qquad \dom \qf h^\rho := \Sobx \rho {X_0, \HS},
\end{equation}
i.e., functions in the domain of $\qlaplacian[\rho]_{X_0}$ fulfill the
usual (now vector-valued) continuity and Kirchhoff
conditions~\eqref{eq:cont} and~\eqref{eq:kirch} on all \emph{inner}
vertices $v \in \ring V_0= \ring Y \cap V_0$. Similarly, we define the
Dirichlet and Kirchhoff $\HS$-valued operators $\qlaplacian[\bd V]_{Y}
\otimes \1$ and $\qlaplacian_{Y} \otimes \1$ via their
quadratic forms defined similarly as in~\eqref{eq:def.qf}, but with
domains $\Sobx {\bd V} {Y} \otimes \HS$ and $\Sob {Y} \otimes \HS$,
respectively. Note that these operators are \emph{decoupled} in the
following sense: Assume that $\HS$ is $r$-dimensional and
\begin{equation*}
  (f_1,\dots, f_r) \cong 
  f \in \Sob {Y} \otimes \HS \cong \underbrace{\Sob{Y}
      \oplus \dots \oplus \Sob{Y}}_{\text{$r$-times}},
\end{equation*}
then $\qlaplacian_Y \otimes \1$ is unitarily equivalent to the
direct sum of $r$ copies of $\qlaplacian_{Y}$, and therefore the
different components \emph{decouple}. The same statement holds for the
Dirichlet Laplacian $\qlaplacian[\bd V]_Y \otimes \1$.

Our crucial observation, already made
in~\cite{lledo-post:07,lledo-post:08} is the following inclusion of
quadratic form domains
\begin{equation*}
  \Sobx {\bd V} {Y} \otimes \HS \subset
  \Sobx \rho {X_0,\HS} \subset
  \Sob {Y} \otimes \HS,
\end{equation*}
implying first, that if $\HS$ is finite-dimensional, then
$\qlaplacian[\rho]_{X_0}$ has purely discrete spectrum (denoted by
$\lambda_k^\rho$, written in ascending order and repeated with respect
to multiplicity). Moreover we have the following assertion proven via
the min-max characterisation of the eigenvalues as
in~\cite{lledo-post:07,lledo-post:08}:
\begin{proposition}
  \label{prp:sp.incl}
  Assume that $\rho$ is a $r$-dimensional representation (i.e., $\dim
  \HS=r$). Then
  \begin{equation*}
    \lambda_k^{\bd V} \ge \lambda_j^\rho \ge \lambda_k, \qquad
    j=(k-1) r+1, \dots, k r.
  \end{equation*}
  In other words, the $j$-th $\rho$-equivariant eigenvalue is enclosed
  in the $k$-th metric Kirchhoff-Dirichlet interval
  \begin{equation*}
    \lambda_j^\rho \in I_k=I_k(Y,\bd V),  \qquad
    j=(k-1) r+1, \dots, k r.
  \end{equation*}
  Moreover, in the equilateral case and for those indices $k$ of case
  BB described in \Prp{kd-int.mg}, the $\rho$-equivariant eigenvalues
  are \emph{independent} of $\rho$, and given by
  $\lambda_j^\rho=(n+1)^2\pi^2$. The corresponding eigenfunctions are
  precisely the topological eigenfunctions of the graph $H$ with
  boundary $\bd V$ and supported in the \emph{interior} of the
  fundamental domain.
\end{proposition}

\subsection{Equivariant discrete Laplacians}
For simplicity, we assume that our discrete graphs have the standard
weights. Let $G=(V,E,\bd) \to G_0=(V_0,E_0,\bd_0)$ be a covering of
discrete graphs with covering group $\Gamma$ and finite quotient graph
$G_0 = G/\Gamma$. Let $\rho$ be a unitary representation of $\Gamma$
on the Hilbert space $\HS$. Denote by
\begin{equation*}
  \lsqr[\rho] {V_0,\HS} :=
   \bigset{\map F V \HS} {F(\gamma \cdot v) = \rho(\gamma) F(v), \quad
     v \in V}
\end{equation*}
the space of $\rho$-equivariant functions.  Again, functions in
$\lsqr[\rho] {V_0,\HS}$ are determined by their values on the vertices
of the quotient $V_0$ (as the notation already indicates). We denote
by $\dlaplacian [\rho]_{G_0}$ the \emph{$\rho$-equivariant} or
\emph{$\rho$-twisted Laplacian} defined as the restriction of
$\dlaplacian_G \otimes \1$ from $\lsqr V \otimes \HS$ onto
$\lsqr[\rho] {V_0,\HS}$.

Let $Y$ be a fundamental domain of the associated metric graph, such
that~\eqref{eq:fd.bd} holds. Now, $Y$ defines a boundary $\bd V$,
which we will also consider as boundary of the discrete graph. Of
course, $\bd V$ depends on the choice of fundamental domain.  Denote
the discrete graph associated to $Y$ by $H$.

As in \Sec{sp.rel}, we denote by $\discr{\mc N}^\rho(\eta) := \ker
(\dlaplacian[\rho]_G - \eta)$ and $\mc N^\rho (\lambda) := \ker
(\qlaplacian[\rho]_X - \lambda)$ the eigenspaces of the equivariant
discrete and metric Laplacian, respectively. Moreover, $\mc
N^\rho_0(\lambda)$ denotes the subspace of $\mc N^\rho(\lambda)$ of
eigenfunctions vanishing at all vertices (see \Def{ef.top}).

For equilateral metric graphs, we have an analogue of
\Prps{sp.rel}{dec.spaces}. Denote by $1$ the trivial representation on
$\HS=\C$ and by $R_{\mathrm a}$ the set of non-trivial involutive
unitary representations on $\HS=\C$, i.e.,
$\rho(\gamma)^{-1}=\rho(\gamma)^*=\rho(\gamma)$ for $\gamma \in
\Gamma$ and $\rho \ne 1$. We also call $R_{\mathrm a}$ the set of
\emph{antisymmetric} representations of $\Gamma$. Note that $R_{\mathrm
  a}$ may be empty.
\begin{proposition}
  \label{prp:sp.rel.eq}
  Assume that the metric covering graph $X \to X_0$ is equilateral
  such that the quotient $X_0$ is compact and connected.
  \begin{enumerate}
  \item Assume that $\lambda \notin \Sigma^\Dir$, then $\map
    {\Phi_\lambda} {\discr{\mc N}^\rho(\mu(\lambda))} {\mc N^\rho
      (\lambda)}$ with $F \mapsto f=\Phi_\lambda F$ as defined in
    \Prp{sp.rel} is an isomorphism. In particular,
    \begin{equation*}
      \lambda \in \spec {\qlaplacian[\rho]_{X_0}}
              \qquad \text{iff} \qquad
      \mu(\lambda) \in \spec {\dlaplacian[\rho]_{G_0}},
    \end{equation*}
    preserving multiplicity.
  \item If $\lambda_n=n^2\pi^2 \in \Sigma^\Dir$, we have the following
    cases:
    \begin{enumerate}
    \item
      \label{case.a}
      If $\rho \notin R_{\mathrm a} \cup \{1\}$ is irreducible then
      ${\mc N}^\rho(\lambda_n) = {\mc N}_0^\rho(\lambda_n)$.
    \item
      \label{case.b}
      If $\rho=1$ is the trivial representation on $\HS=\C$, then
      \begin{equation*}
        {\mc N}^\rho(\lambda_n) = 
        \begin{cases}
          {\mc N}^\rho_0(\lambda_n) & 
                    \text{$n$ odd and $G_0$ not bipartite,}\\
          {\mc N}^\rho_0(\lambda_n) \oplus \C \phi_n & \text{otherwise.}
        \end{cases}
      \end{equation*}
      Here, $\phi_n$ is associated to the graph $G_0$ (see
      \Def{ef.vx0}).
    \item
      \label{case:c}
      If $\rho \in R_{\mathrm a}$ is antisymmetric, then ${\mc
        N}^\rho(\lambda_n) = {\mc N}^\rho_0(\lambda_n) \oplus \C
      \phi_n$ provided $n$ is odd and $G$ has a bipartite fundamental
      domain $H$ such that the connecting vertices $\gamma \cdot v, v
      \in \bd H$ are joined by a path of odd length. In all other
      cases, ${\mc N}^\rho(\lambda_n) = {\mc N}^\rho_0(\lambda_n)$.
      Here, $\phi_n$ is associated to the the bipartite eigenfunction
      of $H$.
    \end{enumerate}
  \end{enumerate}
\end{proposition}
Note that if $G_0$ is bipartite then any fundamental domain $H$ is,
but not vice versa.
\begin{proof}
  The first statement is analogue to the one of \Prp{sp.rel} and can
  be shown similarly as e.g.~in~\cite{post:pre07c}. The proof of the
  second statement is similar to the proofs of \Lem{two.cases} and
  \Prp{dec.spaces}. We only sketch the ideas here. Let $f \in \mc
  N^\rho(\lambda_n)$ be an eigenfunction, interpreted as function on a
  fundamental domain $Y$. Fix $v \in \bd Y$ and let $\gamma \in
  \Gamma$ such that $\gamma \cdot v \in \bd Y$. Note that the set
  $\Gamma_0$ of all such $\gamma$'s generate the group $\Gamma$ (see
  \cite{ratcliffe:94}). Let $p_\gamma$ be a path from $v$ to $\gamma
  \cdot v$ without passing a vertex twice. Denote by $s(p_\gamma)$ the
  number of edges of $p_\gamma$. Then $f(\gamma \cdot v)=(-1)^{n
    s(p_\gamma)} f(v)$. If the fundamental domain is bipartite, then
  $s(\gamma):= s(p_\gamma)$ is independent of the path joining $v$ and
  $\gamma \cdot v$. Note that $s(\gamma)$ may still depend on $v \in
  \bd V$. For $\gamma' \in \Gamma_0 \setminus \{\gamma\}$ we set
  $s(\gamma')=0$.  Now, $\rho_n(\gamma'):= (-1)^{n s(\gamma')}$
  extends to a unitary representation of $\Gamma$ on $\C$.

  The equivariance condition implies that
  \begin{equation*}
    f(v) \in \bigcap_{\gamma_0 \in \Gamma_0} 
                \ker \bigl(\wt \rho_n(\gamma_0) - \id_{\HS}\bigr)
  \end{equation*}
  where $\wt \rho_n (\gamma'):= \rho_n(\gamma') \rho(\gamma)$ is a
  representation on $\HS$. Since $\rho$ is irreducible, $\wt \rho_n$
  is also irreducible. Moreover, if $\rho \ne \rho_n$, then $f(v)=0$
  and vanishes therefore on all vertices, since an irreducible
  representation not in $R_{\mathcal a} \cup \{1\}$ cannot have a
  common eigenvector.  This covers Case~\eqref{case.a}. Otherwise,
  $\rho = \rho_n$ and $\HS=\C$, and in particular, $\rho=1$ if $n$ is
  even or $\rho \in R_{\mathrm a}$ if $n$ is odd.  The other cases
  follow step by step.  Note that in Case~\eqref{case.b}, $n$ odd, it
  follows from the bipartiteness of the quotient graph $G_0$, that $s$
  is even, and therefore $f(v)=\phi_n(v)$ is a vertex-based solution.
\end{proof}

\begin{remark}
  \label{rem:twisted.hom}
  We can equivalently define the Laplacian as $\dlaplacian[\rho]_{G_0}
  = d_\rho^* d_\rho$ where $d_\rho$ is a ``twisted'' exterior
  derivative, defined via
  \begin{equation*}
    \map {d_\rho}{\lsqr[\rho]{V_0,\HS}} {\lsqr E \otimes \HS}, \qquad
    (d_\rho F)_e = F(\bd_+e) - F(\bd_-e).
  \end{equation*}
  Moreover, one can show that the mapping
  \begin{equation*}
    \map{\wt \Psi_n} {\ker d_\rho^*} {{\mc N}_0^\rho(\lambda_n)},
        \qquad
  \end{equation*}
  given by $\wt \Psi_n \eta = f_\eta$, and $f_{\eta,e}(x)=\eta_e
  \sin(n \pi x)$ for $n$ even is an isomorphism, i.e., the topological
  eigenfunctions of the twisted metric graph are related to the
  twisted cohomology $H^1_\rho(X_0,\HS) := \ker d_\rho^*$. Defining
  the corresponding twisted homologies $H_1^\rho(X_0,\HS)$ as in
  \Sec{hom}, we obtain the statement analogue to \Prp{hom.ef}.

  Similarly, for $n$ odd we obtain the corresponding statements for
  the unoriented version $\unor d_\rho$ and the related
  (co-)homologies. We skip the details here, as well as an analysis of
  the twisted Betti numbers, since we do not need the precise spectral
  information of $\dlaplacian[\rho]_{X_0}$ for the existence of gaps.
\end{remark}

We can now carry over the results of \Prp{sp.incl} to discrete graphs.
Note that $\qlaplacian[\bd V]_H$ is equivalent to a square matrix of
size $\abs V-\abs{\bd V}$ where $V=V(H)$. Similarly, $\dlaplacian_H$
is described by an $\abs V \times \abs V$-matrix and
$\dlaplacian[\rho]_{G_0}$ by an matrix of size $r\abs{V_0}$ where
$r=\dim \HS$ and $V_0=V(G_0)$. Moreover, $\abs V -\abs{\bd V} \le
\abs{V_0} \le \abs V$.
\begin{proposition}
  \label{prp:sp.incl.dg}
  Assume that $\rho$ is an $r$-dimensional representation (i.e., $\dim
  \HS=r$). Then\footnote{If $k > \abs V - \abs{\bd V}$, then there are
    no Dirichlet eigenvalues left. In this case, the inequality is
    understood as if we would have set $\mu_k^{\bd V}=2$. This is
    consistent with the definition of the discrete KD intervals (see
    \Def{kd.int.dg}).}
  \begin{equation*}
    \mu_k^{\bd V} \ge \mu_j^\rho \ge \mu_k, \qquad
    j=(k-1) r+1, \dots, k r, \quad
    k=1, \dots, \abs{V_0}.
  \end{equation*}
  In other words, the $j$-th $\rho$-equivariant eigenvalue is enclosed
  in the discrete Kirchhoff-Dirichlet interval
  \begin{equation*}
    \mu_j^\rho \in J_k=J_k(H,\bd V),  \qquad
    j=(k-1) r+1, \dots, k r, \quad
    k=1, \dots, \abs{V_0}.
  \end{equation*}
\end{proposition}
Note that the discrete KD intervals are defined for $k \in
\{1,\dots,\abs V \}$ (see \Def{kd.int.dg}), whereas the
$\rho$-equivariant eigenvalues are given only for $k \le \abs{V_0}$.

%
\section{Residually finite coverings}
\label{sec:cov.gr}
%

We consider now \emph{infinite} coverings with compact quotient graph
and covering group $\Gamma$.

\subsection{Abelian groups}
Let us start with Abelian covering groups $\Gamma$, for which we have
the powerful tool of Floquet-(Bloch)-decomposition. We state the results only for the metric case, the discrete case can be treated similarly. The direct integral decomposition is of the form
\begin{equation*}
  \Lsqr X \cong \int_{\wh \Gamma}^\oplus \Lsqr Y, \qquad
  \qlaplacian_X \cong \int_{\wh \Gamma}^\oplus \qlaplacian[\rho]_{X_0}.
\end{equation*}
Since $\Gamma$ is Abelian, $\rho$ can be parametrised by $\vartheta
\in \R^r$ via $\rho(\gamma)=\e^{\im \vartheta \cdot \gamma}$. We also
write $\lambda_k^\vartheta$ for $\lambda_k^\rho$. For details we refer
to~\cite[Sec.~6]{sunada:pre08} or~\cite{lledo-post:07} and the
references therein.  Moreover, from the direct integral decomposition
and the continuous dependence of $\lambda_k^\rho$ on $\rho$, 
we deduce for the spectrum of the Kirchhoff Laplacian
\begin{equation}
  \label{eq:band}
  \spec {\qlaplacian_X}
  = \bigcup_{\rho \in \wh \Gamma} \spec{\qlaplacian[\rho]_{X_0}}
  = \bigcup_{k \in \N} B_k 
    \qquad \text{where} \qquad
  B_k := \set{\lambda_k^\rho}{\rho \in \wh \Gamma}
\end{equation}
is called the \emph{$k$-th band} and $\lambda_k^\rho$ denotes the
$k$-th eigenvalue of the equivariant Laplacian
$\qlaplacian[\rho]_{X_0}$. The next proposition is a direct
consequence of \Prp{sp.incl}:
\begin{proposition}
  Denote by $I_k=I_k(H,\bd V)$ the metric KD interval of the
  fundamental domain $H$ with vertices $V=V(H)$ and edges $E=E(H)$.
  Then we have
  \begin{equation*}
    \spec {\qlaplacian_X} = \bigcup_{k \in \N} B_k 
    \subset \bigcup_{k \in \N} I_k.
  \end{equation*}
  In particular, the bands $B_k$ with index $k=n \abs E + \abs V +
  1-\alpha_n, \dots, (n+1)\abs E$, ($\alpha_n=1$ if $G$ is bipartite
  or $G$ is not bipartite and $n$ odd, $\alpha_n=0$ otherwise) are
  reduced to points $\{(n+1)^2\pi^2\}$. Moreover, if $\chi(H)=\abs V
  -\abs E \le \alpha_n-1$ then $(n+1)^2\pi^2$ is an eigenvalue of
  infinite multiplicity for $\qlaplacian_X$. The corresponding
  eigenspaces are generated by compactly supported edge-based
  (topological) eigenfunctions of the fundamental domain $H$ and its
  translates.
\end{proposition}
 
\subsection{Residually finite groups}

The following construction of covering graphs is valid for the
discrete and metric case by assuming that the projection respects the
corresponding structure, i.e., they are graph morphisms respecting
orientation in both cases, and additionally, they preserve the length
functions.

Assume that $X_0$ is compact (i.e.  finite for \emph{discrete}
graphs). Moreover, suppose that $\map \pi X {X_0}$ is a covering with
covering group $\Gamma=\Gamma_0$. Corresponding to a normal subgroup
$\Gamma_i \lhd \Gamma$ we associate a covering $\map {\pi_i} X {X_i}$
such that
\begin{equation}
  \label{eq:sub.cov}
  \begin{diagram}
  &                                 & X &                         &\\
  &  \ldTo(2,2)^{\pi_i}_{\Gamma_i}  &   & \rdTo(2,2)^\pi_\Gamma   & \\
  X_i  &                            & \rTo^{p_i}_{\Gamma/\Gamma_i} 
                                                                & & X_0
  \end{diagram}
\end{equation}
is a commutative diagram. The groups under the arrows denote the
corresponding covering groups.

\begin{definition}
  \label{def:res.fin}
  A (countable, infinite) discrete group $\Gamma$ is residually finite
  if there exists a monotonely decreasing sequence of normal subgroups
  $\Gamma_i \lhd \Gamma$ such that
  \begin{equation}
    \label{eq:sub.groups}
    \Gamma=\Gamma_0 \rhd \Gamma_1 \rhd \dots \rhd \Gamma_i \rhd \cdots,  \quad
    \bigcap_{i \in \N} \Gamma_i = \{e\}  \quad \text{and} \quad
    \text{$\Gamma/\Gamma_i$ is finite.}
  \end{equation}
\end{definition}
Suppose now that $\Gamma$ is residually finite. Then there exists a
corresponding sequence of coverings $\map {\pi_i} X {X_i}$ such that
$\map {p_i} {X_i} {X_0}$ is a \emph{finite} covering
(cf.~Diagram~\eqref{eq:sub.cov}). Such a sequence of covering maps is
also called \emph{tower of coverings}.

For more details on residually finite groups we refer
to~\cite{lledo-post:08} and the references therein. The next
proposition is provided by Adachi~\cite{adachi:95} (see
also~\cite[Sec.~5]{lledo-post:08}). We just mention the geometric
meaning of this algebraic condition: The covering space $X$ with
residually finite group can be ``exhausted'' by the finite covering
spaces $X_i$ as one uses in the next proposition.  Its proof can be
redone literally as in the manifold case.
\begin{proposition}
  \label{prp:res.fin}
  Suppose $\Gamma$ is residually finite with the associated sequence
  of coverings $\map {\pi_i} X {X_i}$ and $\map {p_i} {X_i} {X_0}$ as
  in~\eqref{eq:sub.cov}.  Then
  \begin{equation*}
    \spec {\qlaplacian_X} \subseteq 
    \clo {\bigcup_{i \in \N} \spec {\qlaplacian_{X_i}}},
  \end{equation*}
  and the Laplacian $\qlaplacian_{X_i}$ w.r.t.\ the finite covering
  $\map {p_i} {X_i} {X_0}$ has discrete spectrum. Equality holds iff
  $\Gamma$ is amenable.
\end{proposition}

Next we analyse the spectrum of the finite covering $X_i \to X_0$ as
in~\cite{lledo-post:08}.  Note that a fundamental domain for $X \to
X_0$ can also be viewed as fundamental domain for \emph{each} finite
covering $X_i \to X_0$, $i \in \N$.
\begin{proposition}
  \label{prp:fin.group}
  We have
  \begin{equation*}
    \spec {\qlaplacian_{X_i}} = 
       \bigcup_{[\rho] \in \widehat{\Gamma/\Gamma_i}}
                       \spec {\qlaplacian[\rho]_{X_0}},
  \end{equation*}
  where $\qlaplacian[\rho]_{X_0}$ is the equivariant Laplacian
  introduced in \Def{eq.lapl} and $\Gamma/\Gamma_i$ is a finite group
  and $\widehat{G_i}$ its dual, i.e., the set of equivalence classes
  $[\rho]$ of unitary, irreducible representations $\rho$ of $\Gamma$.
\end{proposition}
In particular, we have:
\begin{theorem}
  \label{thm:brack.mg} 
  Suppose $X \to X_0$ is a $\Gamma$-covering of (not necessarily
  equilateral) metric graphs with fundamental domain $Y$, where
  $\Gamma$ is a residually finite group.  Then
  \begin{equation*}
    \spec {\qlaplacian_X} \subset \bigcup_{k \in \N} I_k =: I,
  \end{equation*}
  where $I_k :=[\lambda_k, \lambda_k^{\bd V}]$ are the
  Kirchhoff-Dirichlet intervals associated to the fundamental domain
  $Y$ with boundary vertices $\bd V$ (see \Def{kd.int.mg}). In
  particular, if $M$ is an interval with $M \cap I = \emptyset$, then
  $M \cap \spec {\qlaplacian_X} = \emptyset$ (i.e., $M$ is a spectral
  gap).

  Moreover if $G$ is bipartite, we have the spectral inclusion
  \begin{equation*}
    \spec {\qlaplacian_X} \subset \hat I
         \qquad \text{where} \qquad
    \hat I := \bigcup_{n=0}^\infty \tau_n(I \cap K_n) \cap (I \cap K_n)
  \end{equation*}
  and where $\tau_n$ is the spectral symmetry defined in \Prp{bip.mg}.
\end{theorem}
\begin{proof}
  We have
  \begin{equation*}
    \spec {\qlaplacian_X} \subseteq 
    \clo {\bigcup_{i \in \N} \spec {\qlaplacian_{X_i}}} =
    \clo {\bigcup_{i \in \N}
          \bigcup_{[\rho] \in \wh {\Gamma_i}} 
                \spec {\qlaplacian[\rho]_{X_0}}}
    \subseteq  
    \overline {\bigcup_{k \in \N} I_k } = \bigcup_{k \in \N} I_k,
  \end{equation*}
  where we used \Prpss{sp.incl}{res.fin}{fin.group}. The results for
  $\hat I$ follow from the spectral symmetry for $\qlaplacian_X$.
\end{proof}

Similarly, in the discrete case, we conclude from
\Prpss{sp.incl.dg}{res.fin}{fin.group}:
\begin{theorem}
  \label{thm:brack.dg} 
  Suppose $G \to G_0$ is a $\Gamma$-covering with fundamental
  domain $H$, where $\Gamma$ is a residually finite group, then
  \begin{equation*}
    \spec {\dlaplacian_G} \subset \bigcup_{k=1}^{\abs V} J_k =: J,
  \end{equation*}
  where $J_k :=[\mu_k, \mu_k^{\bd V}]$ are the discrete Kirchhoff-Dirichlet
  intervals associated to the fundamental domain $H$ with boundary vertices
  $\bd V$ (see \Def{kd.int.dg}).  In particular, if $M \cap J = \emptyset$,
  then $M \cap \spec {\dlaplacian_G} = \emptyset$ (i.e., $M$ is a spectral
  gap).

  Moreover if $G$ is bipartite, we have the spectral inclusion
  \begin{equation*}
    \spec {\dlaplacian_G} \subset \hat J
         \qquad \text{where} \qquad
    \hat J := \theta(J) \cap J
  \end{equation*}
  and where $\theta(\mu)=2-\mu$ is the spectral symmetry defined in
  \Prp{bipartite}.
\end{theorem}
We refer to $I$ and $J$ as the \emph{KD spectrum} and to $\hat I$ and
$\hat J$ as the \emph{symmetrised} KD spectrum.

Let us mention separately the case when $\Gamma$ is amenable: 
\begin{theorem}
  \label{thm:amenable}
  Assume that the covering group of the covering is amenable. Then the
  number of components of $\spec{\qlaplacian_X}$ resp.\ of
  $\spec{\dlaplacian_G}$ is at least as large as the number of
  components of $I$ resp.\ $J$ or $\hat I$ resp\ $\hat J$ in the
  bipartite case.
\end{theorem}
\begin{proof}
  The assertion follows from the fact that, due to amenability, we
  have equality in~\Prp{res.fin}. In particular, the spectrum of
  $\qlaplacian_{X_0}$ is contained in $\qlaplacian_X$. Moreover, the
  quotient spectrum is just the spectrum of $\qlaplacian[\rho]_{X_0}$
  with the trivial representation $\rho=1$. Therefore, the $k$-th
  eigenvalue $\lambda_k(X_0)$ of $\qlaplacian_{X_0}$ is contained in
  $\spec {\qlaplacian_X}$ and also in the $k$-th KD interval $I_k$ due
  to \Prp{sp.incl}. The discrete case follows similarly.
\end{proof}
%
\section{Examples: covering graphs with spectral gaps}
\label{sec:ex}
%

In this section we present several examples for which the KD intervals
already guarantee the existence of spectral graphs. In some cases, the
symmetrised KD spectrum is even \emph{equal} to the $\Z^r$-periodic
spectrum, see the bipartite examples below. For the concrete examples
one only needs to calculate the spectra of the matrices associated to
the discrete operators (see \Eq{mat.std}) on a suitable chosen
fundamental domain. For brevity, we skip the corresponding spectral
results for metric graphs, since they can be obtained straightforward
by the results of the previous sections.

In \ExS{bip}{self-loops}, we consider ``small'' $\Z$-periodic graphs
in order to show how our method works in simple examples in which the
periodic spectrum can also be calculated directly. One can see that
the KD intervals give ``good'' estimates of the actual location of the
bands only for the first and second band.  For a larger number of gaps
guaranteed by the KD intervals, one should consider graphs with a
smaller ratio $\abs{\bd V}/\abs V$.  Of course, our method is more
interesting for non-abelian (residually finite) groups with more than
one generator, see \Ex{2gen}.

Note that the choice of fundamental domain is arbitrary, and that the
definition of the KD intervals will (in general) depend on the choice
of the fundamental domain. Therefore, it might happen, that a ``good''
choice of fundamental domain leads to a union of the KD intervals
having gaps. We do not precise the meaning of ``good'' here, but as in
the case of manifolds and Schr\"odinger operators (see
e.g.~\cite{hempel-post:03,lledo-post:07,lledo-post:08}) the
fundamental domain should have ``small'' boundary in order to decouple
from its neighbours.  In our context, this means that a fundamental
domain $H$ should contain a large number of vertices $V(H)$ and edges
$E(H)$ compared to the number of boundary vertices $\bd V$, see
\Exs{mult-edge}{self-loops}.

We start with a bipartite example already used in \Fig{sp.rel}.
\begin{example}
  \label{ex:bip}
  Let $G \to G_0$ be the periodic graph with fundamental domain $H$ as
  given in \Fig{gr.bip}. The spectrum of the discrete (Dirichlet)
  Laplacian is
  \begin{equation*}
    \spec{\dlaplacian_H} = \{0,1,1,1,2\} \Und
    \spec{\dlaplacian[\bd V]_H} 
     = \bigg\{1-\frac 1 {\sqrt 3}, 1,1+\frac 1 {\sqrt 3} \bigg\},
  \end{equation*}
  resp., where repeated numbers correspond to multiple
  eigenvalues, so that the KD intervals are
  \begin{equation*}
    J_1=\bigg[0, 1-\frac 1 {\sqrt 3}\bigg], \quad
    J_2=\{1\}, \quad
    J_3=\bigg[1, 1+\frac 1 {\sqrt 3}\bigg], \quad
    J_4=[1,2].
  \end{equation*}
  The equivariant spectrum for $\rho(\gamma)=\e^{\im \vartheta \gamma}$ is
  \begin{equation*}
    \spec{\dlaplacian[\vartheta]_{G_0}}
     = \bigg\{ 1- \sqrt{\frac {2 + \cos \vartheta} 3},
            1,1, 1 + \sqrt{\frac {2 + \cos \vartheta} 3} \Bigg\}.
  \end{equation*}
  In particular, the bands $\discr B_k := \set{\lambda^\vartheta_k}{\vartheta
    \in [0,2\pi]}$  (see \Eq{band}) are
  \begin{equation*}
    \discr B_1=\bigg[0, 1-\frac 1 {\sqrt 3}\bigg], \quad
    \discr B_2= \discr B_3=\{1\}, \quad
    \discr B_4=\bigg[1+\frac 1 {\sqrt 3},2 \bigg].
  \end{equation*}
  We see that the first and second band agree with the corresponding
  KD intervals. In particular, the KD intervals detect the first gap
  $(1-1/\sqrt 3,1)$ \emph{precisely}. Moreover, the second KD interval
  reduces to a point as well as the second band. But the third KD
  interval is too rough, and the second gap is not detected. See also
  \Fig{sp.rel} for the spectral relation with the corresponding metric
  graph. Nevertheless, since $G$ is bipartite, we also have the
  spectral inclusion for the symmetrised KD spectrum $\hat J$. Here,
  we even have the equality $\spec {\dlaplacian_G} = \hat J$ (see
  \Thm{brack.dg}), showing that the KD intervals can give the actual
  spectrum of the covering using the spectral symmetry for bipartite
  graphs.
\begin{figure}[h]
  \centering
\begin{picture}(0,0)%
  \includegraphics{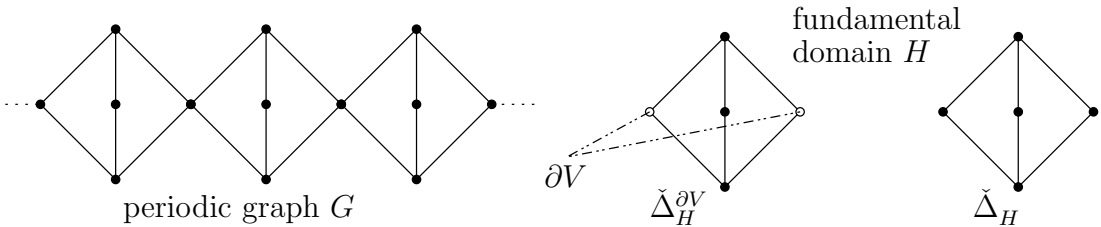}%
\end{picture}%
\setlength{\unitlength}{4144sp}%
\begin{picture}(6570,1294)(214,-650)
\put(946,-601){periodic graph $G$}%

\put(4951,329){domain $H$}%

\put(4951,524){fundamental}%

\put(3466,-421){$\bd V$}%

\put(4096,-601){$\dlaplacian[\bd V]_H$}%

\put(6031,-601){$\dlaplacian_H$}%

\end{picture}%
  
  \caption{A bipartite graph. The related spectral information of this
    graph is visualised in \Fig{sp.rel}. The full vertices correspond
    to Kirchhoff conditions and the open vertices correspond to
    Dirichlet conditions for the associated metric graph.}
  \label{fig:gr.bip}
\end{figure}
\end{example}

The next example is a non-bipartite one:
\begin{example}
  \label{ex:non-bip}
  Let $G \to G_0$ be the periodic graph with fundamental domain $H$ as
  given in \Fig{gr.non-bip}. The spectrum of the discrete (Dirichlet)
  Laplacian is
  \begin{equation*}
    \spec{\dlaplacian_H} 
     = \bigg\{0,
      \frac {7-\sqrt {13}}6,
      \frac 43, \frac 43,
      \frac {7+\sqrt {13}}6 \bigg\} \Und
    \spec{\dlaplacian[\bd V]_H} 
     = \bigg\{\frac 13, \frac 43, \frac 43 \bigg\},
  \end{equation*}
  resp., where repeated numbers correspond to multiple eigenvalues, so
  that the KD intervals are
  \begin{equation*}
    J_1=\bigg[0, \frac 13\bigg], \quad
    J_2=\bigg[\frac {7-\sqrt {13}}6, \frac 43\bigg], \quad
    J_3=\bigg\{ \frac 43 \bigg\}, \quad
    J_4=\bigg[\frac 43, 2\bigg].
  \end{equation*}
  The spectrum of the periodic operator is given by the bands
  \begin{equation*}
    \discr B_1= \bigg[0, 1-\frac{\sqrt 5} 3\bigg], \quad
    \discr B_2= \bigg[\frac 23, \frac 43 \bigg], \quad
    \discr B_3= \bigg\{ \frac 43 \bigg\}, \quad
    \discr B_4= \bigg[\frac 43, 1+\frac{\sqrt 5} 3 \bigg].
  \end{equation*}
  Here, only the degenerated band $\discr B_3$ agrees with the KD
  interval.  The corresponding eigenfunction is indicated by the
  values in \Fig{gr.non-bip}. Nevertheless, the maximal spectral gap
  in this example $(1-\sqrt 5/3, 2/3) \approx(0.25,0.66)$ is detected
  approximately by the KD interval giving the spectral gap $(1/3,
  (7-\sqrt{13})/6) \approx (0.33, 0,57)$. The fourth KD interval gives
  a too rough upper bound, see also \Fig{sp.rel.n-bp} for the spectral
  relation with the corresponding metric graph.

  Note that the (metric) KD intervals do not detect the gap between
  the nineth and tenth band, the KD intervals even overlap (see
  \Fig{sp.rel.n-bp}), whereas the gap between the sixth and seventh
  band is recognised.
\begin{figure}[h]
  \centering
\begin{picture}(0,0)%
  \includegraphics{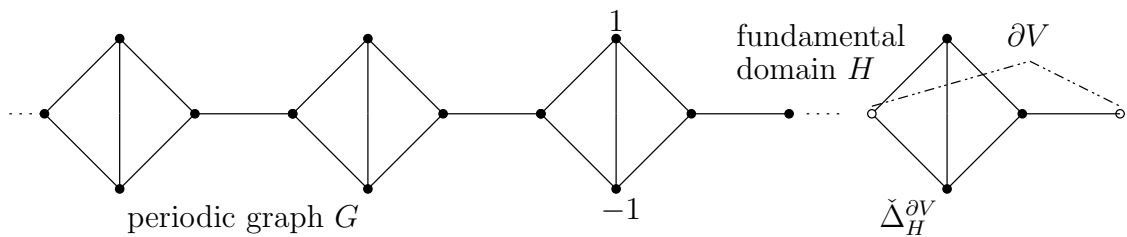}%
\end{picture}%
\setlength{\unitlength}{4144sp}%
\begin{picture}(6705,1327)(214,-650)
\put(946,-601){periodic graph $G$}%

\put(4591,284){domain $H$}%

\put(4591,479){fundamental}%

\put(6211,479){$\bd V$}%

\put(5446,-601){$\dlaplacian[\bd V]_H$}%

\put(3781,-556){$-1$}%

\put(3826,569){$1$}%

\end{picture}%
  \caption{A non-bipartite graph. The related spectral information of
    this graph is visualised in \Fig{sp.rel.n-bp}. The Laplacian on
    $H$ without Dirichlet conditions is not plotted here. The two
    values $\pm1$ at the vertices indicate the eigenfunction
    associated to the eigenvalue $\lambda_3=4/3$, independent of
    $\vartheta$; the other vertex values being $0$.}
  \label{fig:gr.non-bip}
\end{figure}
\end{example}

The following example (see also~\cite{ael:94}, where the band-gap
ratio of such ``onion-like'' periodic metric graphs is considered)
gives an idea of how to generate gaps by multiple edges:
\begin{example}
  \label{ex:mult-edge}
  Let $G \to G_0$ be the periodic graph with fundamental domain $H$ as
  given in \Fig{gr.mult-edge} having $r$ repeated edges. The spectrum
  of the discrete (Dirichlet) Laplacian is
  \begin{equation*}
    \spec{\dlaplacian_H} 
     = \bigg\{0,
      1 - \frac 1 {r+1},
      1 + \frac 1 {r+1}, 
      2 \bigg\} \Und
    \spec{\dlaplacian[\bd V]_H} 
     = \bigg\{\frac 1 {r+1}, 2 - \frac 1{r+1} \bigg\},
  \end{equation*}
  respectively.  The KD intervals are
  \begin{equation*}
    J_1=\bigg[0, \frac 1 {r+1} \bigg],  \quad
    J_2=\bigg[1 - \frac 1 {r+1}, 2 - \frac 1{r+1}\bigg]
          \und
    J_3=\bigg[1 + \frac 1 {r+1}, 2 \bigg].
  \end{equation*}
  Note that as far as $r \ge 2$, we have spectral gaps between the
  first and second KD interval. Moreover, the KD intervals reduce to
  the point $\{0\}$ for $k=1$ and to the interval $[1,2]$ for $k=2,3$
  as $r \to \infty$.  The spectrum of the periodic operator is given
  by the bands
  \begin{equation*}
    \discr B_1= \bigg[0, \frac 1 {r+1} \bigg], \quad
    \discr B_2= \bigg[1- \frac 1 {r+1}, 1 + \frac 1 {r+1} \bigg], \quad
    \discr B_3= \bigg[2 - \frac 1 {r+1}, 2 \bigg],
  \end{equation*}
  and only the first KD interval $J_1$ agrees with the first band
  $\discr B_1$. Note that in this case, the periodic and antiperiodic
  equivariant eigenvalues ($\vartheta=0$ and $\vartheta=\pi$) give
  already the band edges. For groups with more than one generator,
  the band edges need not to be on the boundary of the Brillouin zone,
  see~\cite{hksw:07} and appear as KD eigenvalues, but with
  alternating role ($B_k=[\lambda_k^0,\lambda_k^\pi]$ for $k=1,3$ and
  $B_2=[\lambda_2^\pi,\lambda_2,0]$). This phenomena also appears for
  Schr\"odinger operators (see~\cite{kuchment-post:07} and the
  references therein).

  Nevertheless, the graph is bipartite, so we can use the spectral
  symmetry and indeed, we have equality $\spec {\dlaplacian_G} = \hat
  J$. Again, the symmetrised KD spectrum gives already the precise
  spectral information.
  \begin{figure}[h]
    \centering
\begin{picture}(0,0)%
\includegraphics{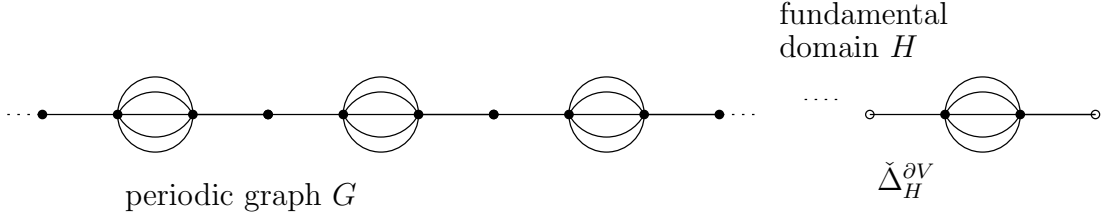}%
\end{picture}%
\setlength{\unitlength}{4144sp}%
\begingroup\makeatletter\ifx\SetFigFont\undefined%
\gdef\SetFigFont#1#2#3#4#5{%
  \reset@font\fontsize{#1}{#2pt}%
  \fontfamily{#3}\fontseries{#4}\fontshape{#5}%
  \selectfont}%
\fi\endgroup%
\begin{picture}(6570,1249)(214,-605)
\put(4861,329){domain $H$}%

\put(4861,524){fundamental}%

\put(946,-556){periodic graph $G$}%

\put(5446,-466){$\dlaplacian[\bd V]_H$}%

\end{picture}%
    \caption{Generating gaps by multiple edges. Here, we replaced the
      middle edge by $r=5$ edges.}
    \label{fig:gr.mult-edge}
  \end{figure}
\end{example}

  A similar result holds by attaching self-loops to a graph:
\begin{example}
  \label{ex:self-loops}
  Let $G \to G_0$ be the periodic graph with fundamental domain $H$
  being a line graph with three vertices and two edges, and $r$ loops
  attached to the middle vertex. The boundary vertices have degree
  $1$, and the middle vertex has degree $2(r+1)$. Note that $G$ is not
  bipartite as long as $r \ge 1$. The spectrum of the discrete
  (Dirichlet) Laplacian can be calculated as
  \begin{equation*}
    \spec{\dlaplacian_H} 
     = \bigg\{0,
      1,
      1 + \frac 1 {r+1}
      \bigg\} \Und
    \spec{\dlaplacian[\bd V]_H} 
     = \bigg\{\frac 1 {r+1} \bigg\}.
  \end{equation*}
  The KD intervals are
  \begin{equation*}
    J_1=\bigg[0, \frac 1 {r+1} \bigg]  \Und
    J_2=[1, 2]
  \end{equation*}
  Note that as far as $r \ge 1$, we have a spectral gap between the
  two KD intervals. Moreover, the first KD intervals reduce to the
  point $\{0\}$ as $r \to \infty$.  The spectrum of the periodic
  operator is given by the bands
  \begin{equation*}
    \discr B_1= \bigg[0, \frac 1 {r+1} \bigg] \Und
    \discr B_2= \bigg[1, 1 + \frac 1 {r+1}\bigg]
  \end{equation*}
  and again, only the first KD interval $J_1$ agrees with the first
  band $\discr B_1$.
\end{example}

We finally present an example with two generators. This example serves
also as an example for coverings with non-abelian groups.
\begin{example}
  \label{ex:2gen}
  Let $G \to G_0$ be the $\Z^2$-periodic graph with fundamental domain
  $H$ as given in \Fig{gr.2gen}. One can also construct other
  coverings associated to a group with two generators by gluing
  together appropriate copies of the fundamental domain according to
  the Cayley graph associated with this generator set.  The discrete
  (Dirichlet) Laplacian is
  \begin{gather*}
    \spec{\dlaplacian_H} 
     = \bigg\{0,
      1-\frac 1 {\sqrt 2},
      1-\frac 1 {\sqrt 2},
      \frac 12,
      1,1,1,1,1,
      \frac 32,
      1+\frac 1 {\sqrt 2},
      1+\frac 1 {\sqrt 2},
      2 \bigg\}, \\
    \spec{\dlaplacian[\bd V]_H} 
     = \bigg\{
      1-\frac {\sqrt 3} 2,
      \frac 12,
      \frac 12,
      1,1,1,
      \frac 32,
      \frac 32,
      1+\frac {\sqrt 3} 2
      \bigg\},
  \end{gather*}
  resp., where repeated numbers correspond to multiple eigenvalues, so
  that the KD intervals are
  \begin{gather*}
    J_1=\bigg[0, 1-\frac {\sqrt 3} 2\bigg] \approx [0,0.13],  \quad
    J_2=J_3=
      \bigg[ 1-\frac 1 {\sqrt 2}, \frac 12 \bigg] \approx [0.29,0.5]\\
    J_4=\bigg[\frac 12,1 \bigg], \quad
    J_5=J_6=\{1\}, \quad
    J_7=J_8=\bigg[1,\frac 32\bigg] \und
    J_9=\bigg[1, 1 +\frac {\sqrt 3} 2\bigg] \approx [1,1.87].  \quad
  \end{gather*}
  It is easily seen that there is a spectral gap only between the first and
  second KD interval. All other intervals overlap. Note that the graph is
  bipartite, so there there is another gap due to the spectral
  symmetry. 

  Here, we can also calculate the spectrum of the $\Z^2$-periodic
  graph using the Floquet theory~\eqref{eq:band}.
  The periodic ($\vartheta=(0,0)$) and antiperiodic
  $(\vartheta=(\pi,\pi)$) spectrum is given by
  \begin{gather*}
    \spec{\dlaplacian[(0,0)]_{G_0}} 
     = \bigg\{0,
      \frac 12,
      \frac 12,
      \frac 12, 1,1,1,
      \frac 32,
      \frac 32,
      \frac 32,
      2 \bigg\}, \\
    \spec{\dlaplacian[(\pi,\pi)]_{G_0}} 
     = \bigg\{
      1-\frac {\sqrt 3} 2,
      1-\frac 1 {\sqrt 2},
      1-\frac 1 {\sqrt 2},
      1,1,1,1,1,
      1+\frac 1 {\sqrt 2},
      1+\frac 1 {\sqrt 2},
      1+\frac {\sqrt 3} 2
      \bigg\},
  \end{gather*}
  respectively. Due to the continuous dependence on $\theta$ and the
  connectedness of $\theta \in [0,2\pi]^2$, we conclude that the
  $k$-th band contains the interval $\wt J_k$ given by the minimum and
  maximum of the $k$-th periodic and antiperiodic eigenvalues. But for
  $k=1,\dots,6$, the interval $\wt J_k$ is already the $k$-th KD
  interval $J_k$, so that $B_k=J_k$ for these $k$.  Moreover, $\discr
  B_5$ and $\discr B_6$ are flat bands. Using the spectral symmetry
  from the bipartiteness, we conclude that the symmetrised KD
  spectrum $\hat J$ already give the spectrum of the $\Z^2$-periodic
  operator $\dlaplacian_G$. In particular, the KD intervals give an
  efficient method to calculate the spectrum of the $\Z^2$-periodic
  Laplacian with a minimum of calculations needed: We only have to
  find ``good'' candidates for $\vartheta$, and do not need the
  spectrum of $\dlaplacian[\vartheta]_{G_0}$ for general $\vartheta$.

  \Thm{brack.dg} assures that $(1-\sqrt 3/2, 1-1\/\sqrt 2)$ and
  $(1+1/\sqrt 2, 1+\sqrt 3/2)$ never belongs to the spectrum of
  \emph{any} covering having $H$ as fundamental domain, in particular
  for the tree-like graph with covering group $\Z * \Z$, the free
  group with two generators. Moreover, $1$ is always an eigenvalue.
\begin{figure}[h]
  \centering
\begin{picture}(0,0)%
\includegraphics{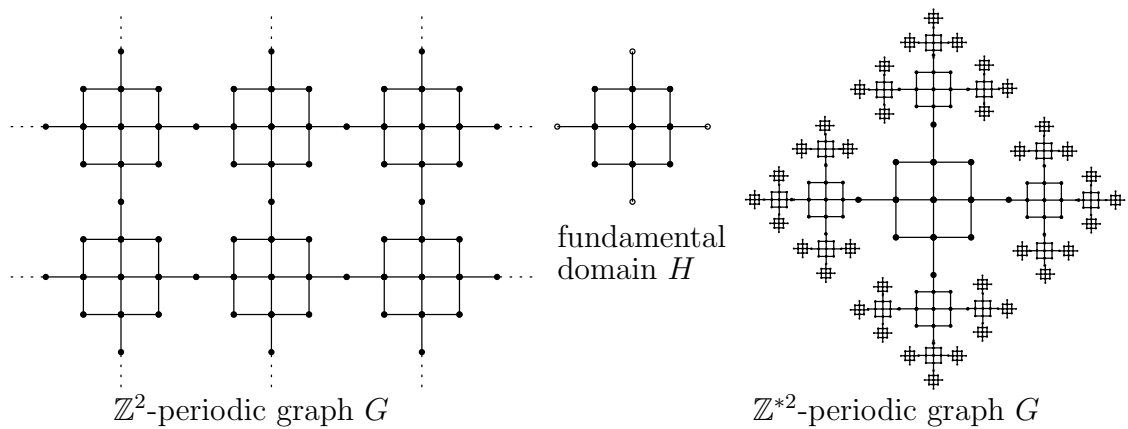}%
\end{picture}%
\setlength{\unitlength}{4144sp}%
\begin{picture}(6702,2518)(79,-1820)
\put(3376,-931){domain $H$}%

\put(3376,-736){fundamental}%

\put(721,-1771){$\Z^2$-periodic graph $G$}%

\put(4546,-1771){$\Z^{*2}$-periodic graph $G$}%

\end{picture}%
  \caption{An example with a covering group having two generators. The
    fundamental domain has 13~vertices and four boundary vertices. On
    the left, the covering graph with Abelian group is plotted. In
    this case, the symmetrised KD spectrum $\hat J$ already give the
    spectrum of the covering operator. On the right, we have a
    $\Gamma$-covering with $\Gamma=\Z^{*2}=\Z * \Z$, the free
    (non-abelian) group with two generators. Here, we only have the
    spectral estimate $\spec {\dlaplacian_{G^*}} \subset \hat J$.}
  \label{fig:gr.2gen}
\end{figure}
\end{example}

\begin{remark}
  \label{rem:kir.opt}
  We could also use (vertex) \emph{Neumann} conditions as lower bound
  on the equivariant metric eigenvalue instead of the \emph{Kirchhoff}
  ones. A function $f$ satisfies the (vertex) Neumann condition in a
  vertex $v \in \bd V$ iff $f_e'(v)=0$ for each edge $e \in E_v$.
  Denote by $\Sobx {\bd V,\Neu} X$ the space of functions $f \in \Sobx
  \max X$ being continuous in each \emph{inner} vertex, i.e., we do
  not assume continuity at boundary vertices. Now, we have the
  additional inclusion $\Sob X \subset \Sobx {\bd V,\Neu} X \subset
  \Sobx \max X$ in~\eqref{eq:sob.incl} and the opposite inequality for
  the eigenvalues, and a similar statement as in \Prp{sp.incl} with
  the Kirchhoff eigenvalue replaced by the vertex Neumann one as lower
  bound. But a direct calculation of the corresponding eigenvalues
  (e.g.~in \Ex{bip}) shows, that the corresponding Neumann-Dirichlet
  intervals do not reveal the spectral gap.

  That our Kirchhoff-Dirichlet bracketing is optimal is shown in
  \Exs{bip}{2gen}, where the (symmetrised) KD spectrum is
  \emph{exactly} the spectrum of the $\Z^2$-periodic graph (and not
  only a superset). 
\end{remark}

\section{Outlook}

We only considered simple examples in which the eigenvalue bracketing
guarantees the existence of spectral gaps. It would be interesting to
provide quantities estimating the actual number of gaps (at least for
Abelian groups $\Gamma=\Z^r$ or amenable groups). As mentioned above,
a naive guess would be that the ratio $\abs{\bd V}/\abs {E(H)}$ is
related to the number of graphs in the union of the KD intervals (the
smaller the ratio is, the more gaps should open up). Moreover,
decorations of the graph (like multiple edges or loops, see
also~\cite{aizenman-schenker:00}) should provide examples with open
gaps, as the examples in \Sec{ex} indicate.  Again, a more systematic
treatment would be interesting.

If the groups $\Gamma$ of the covering is amenable and residually
finite, we provide a lower bound on the number of spectral gaps. The
amenability condition is only needed in order to assure that each KD
interval contains at least one spectral point (namely, an eigenvalue
of the quotient space). This condition might be weakened, but it is a
priori not clear what representation $\rho$ leads to an equivariant
eigenvalue inside the KD interval (see \Prps{res.fin}{fin.group}). In
the case of manifolds, we guaranteed the existence of spectrum inside
the Neumann-Dirichlet intervals for residually finite groups by the
fact that the Dirichlet and Neumann eigenvalues of a suitable chosen
fundamental domain were close to each other
(cf.~\cite[Thm.~3.3]{lledo-post:08}). An upper bound is given once the
covering group has positive Kadison constant (see~\cite{sunada:92}).

Homology groups have also been used for metric graph Laplacians with
\emph{magnetic field}, see~\cite{kostrykin-schrader:03} for details.
The type of spectrum for 
magnetic Laplacians on a metric equilateral square lattice was
analysed in~\cite{bgp:07}, and, in particular, for irrational flux,
the spectrum has Cantor structure.  Magnetic Laplacians may be seen as
a generalisation of equivariant Laplacians for Abelian coverings
treated in detail in \Sec{eq.lapl}.  It would be interesting to see
how the eigenvalue bracketing can be applied to this case in order to
make non-trivial statements about the nature of the spectrum of
discrete and metric magnetic Laplacians.

Another point we do not address here is the appearance of
``degenerated'' bands, i.e., eigenvectors localised inside a
fundamental domain leading to a spectral band reduced to a point. For
metric graphs, this often happens for the exceptional values
$\lambda=n^2 \pi^2$, but this fact can also happen away from these
points, and therefore also for the discrete graph (see
\Exs{non-bip}{2gen}).  Moreover, we do not analyse the band-gap ratio
which may be estimated by from above by the corresponding ratio for the
KD intervals (see~\cite{ael:94} and \Ex{mult-edge}).



\end{document}